\documentclass[11pt,reqno]{amsart}
\allowdisplaybreaks[4]
\usepackage{amssymb}
\usepackage{mathrsfs}
\usepackage{amsthm}
\usepackage{amsmath}
\usepackage{graphicx}
\usepackage{subcaption}
\numberwithin{equation}{section}
\usepackage{bm}
\usepackage[T1]{fontenc}
\usepackage{bbm}
\usepackage{graphicx}
\usepackage{booktabs}
\usepackage{float}
\usepackage{mathrsfs}
\usepackage{geometry}
\usepackage{color}
\usepackage[dvipsnames, svgnames, x11names]{xcolor}
\usepackage{hyperref}
\usepackage{aliascnt}
\usepackage{cleveref}

\theoremstyle{plain}
\newtheorem{lemma}{Lemma}[section]

\newaliascnt{theorem}{lemma}
\newtheorem{theorem}[theorem]{Theorem}
\aliascntresetthe{theorem}
\crefname{theorem}{theorem}{theorems}
\Crefname{theorem}{Theorem}{Theorems}

\newaliascnt{cor}{lemma}

\aliascntresetthe{cor}
\crefname{cor}{corollary}{corollaries}
\Crefname{cor}{Corollary}{Corollaries}

\newaliascnt{prop}{lemma}
\newtheorem{prop}[prop]{Proposition}
\aliascntresetthe{prop}
\crefname{prop}{proposition}{propositions}
\Crefname{prop}{Proposition}{Propositions}

\newtheorem{rem}{Remark}[section]

\Crefname{rem}{Remark}{Remarks}
\Crefname{defn}{Definition}{Definitions}
\Crefname{assp}{Assumption}{Assumptions}

\newcommand{\E}{\mathbb{E}}
\newcommand{\N}{\mathbb{N}}

\newcommand{\PP}{\mathbb{P}}
\newcommand{\RR}{\mathbb{R}}

\newcommand{\op}{\text{op}}

\def\<{\langle}
\def\>{\rangle}

\def\nn{\nonumber}

\def\bc{\begin{center}}       \def\ec{\end{center}}
\def\ba{\begin{array}}        \def\ea{\end{array}}
\def\be{\begin{equation}}     \def\ee{\end{equation}}
\def\bea{\begin{eqnarray}}    \def\eea{\end{eqnarray}}
\def\beaa{\begin{eqnarray*}}  \def\eeaa{\end{eqnarray*}}

\begin{document}
	\title[Almost sure asymptotics of error process]{Almost sure asymptotics of numerical errors for SDEs with additive noise}
	\subjclass[2020]{60F15; 60H35;  65C30}
	\author{Chuchu Chen, Xinyu Chen, Jialin Hong}
	\address{State Key Laboratory of Mathematical Sciences, Academy of Mathematics and Systems Science, Chinese Academy of Sciences, Beijing 100190, China,\and School of Mathematical Sciences, University of Chinese Academy of Sciences, Beijing 100049, China}
	\email{chenchuchu@lsec.cc.ac.cn; chenxinyu@amss.ac.cn; hjl@lsec.cc.ac.cn}
	\thanks{This work is funded by the National Key R\&D Program of China under Grant (No. 2024YFA1015900),  by the National Natural Science Foundation of China (No. 12461160278,  No. 12471386, and No. 12031020), and by Youth Innovation Promotion Association CAS}
	\begin{abstract}
		In this paper,  we characterize the sharp almost sure asymptotics for the numerical error of the Euler--Maruyama method applied to stochastic differential equations with additive noise.
        Specifically, we establish, to our knowledge, the first sharp laws of the iterated logarithm for Euler--Maruyama errors in this setting along geometric mesh refinements $(N_i)_{i\in\N^+}=(\alpha^i)_{i\in\N^+}$ for any integer $\alpha\geq 2$.  These results identify $(N_i)^{-1}\sqrt{\log\log N_i}$ as the sharp almost sure error scale and provide an explicit expression for the path-dependent constants governing the asymptotic error envelopes. The LIL constants also admit a natural interpretation through the asymptotic error distribution, specifically through the conditional covariance matrix of the stable limit of the rescaled error given the driving Wiener process. A central challenge in establishing these results is to characterize the almost sure fluctuations of  Skorohod integrals arising in the error expansion, which are generally anticipating, non-Gaussian, and dependent across refinement levels.  The analysis also provides a framework for studying analogous asymptotics for other numerical methods.
	\end{abstract}

	\keywords{Almost sure asymptotics $\cdot$ Stochastic differential equations $\cdot$ Error analysis $\cdot$ Law of the iterated logarithm}
	\maketitle

	\section{Introduction} 
	Almost sure error analysis provides a natural framework for assessing the accuracy of numerical approximations of individual solution trajectories and a theoretical foundation for reliable pathwise approximation of solutions to stochastic differential equations (SDEs). Such convergence results for numerical methods  have been established in several works; see, e.g.,  \cite{Bichteler1981AOP, Huang2025,  Gyongy1998PA, Kloeden2007LMS, Shardlow2016BIT}. In particular,  \(L^p\)-strong convergence estimates yield almost sure error bounds through the Borel--Cantelli lemma. To illustrate this connection, let $X$ denote the exact solution of an SDE on the time interval $[0,T]$, and $X^N$  denote its numerical approximation with $h=N^{-1}$ for $N\in\N^+$. Suppose that, for some $\tau>0$ and  every $p\geq1$,
	\begin{align*}
		\Big(
		\mathbb E\big[
		|X(T)-X^N(T)|^p
		\big]
		\Big)^{1/p}
		\leq C_{p,T} h^\tau
	\end{align*}
	where $C_{p,T}$ is independent of $N$. Then for every $\varepsilon\in(0,\tau)$,
		\begin{align*}
		|X(T)-X^N(T)|
		\leq C_{\varepsilon,T}(\omega)h^{\tau-\varepsilon},
		\quad\text{a.s.},
	\end{align*} 
	where \(C_{\varepsilon,T}\) is an almost surely finite random variable independent of \(N\). However, these almost sure upper bounds involve an arbitrarily small loss in the convergence order and do not identify the sharp fluctuation scale of the numerical error. This raises a natural question: {\it can one determine the sharp almost sure fluctuation scale as $N\to\infty$, together with the corresponding  constants governing the asymptotic error envelope? }

	In this work, we address this question by characterizing the almost sure asymptotics of the error process along geometric refinement sequences through the law of the iterated logarithm (LIL).  The LIL is a fundamental probabilistic limit theorem that characterizes the sharp asymptotic envelope of pathwise fluctuations by identifying the appropriate scaling function and the corresponding limit constant; see, e.g., \cite{Chenchen2024, Chenchen2025, schilling2014brownian}.
	Specifically,  we consider the $d$-dimensional SDEs with additive noise
	\begin{align} \label{equ_nonliearGeneral}
		dX(t) = f(X(t))dt + \sigma dW(t), \quad X(0)=x_0\in\RR^d,
	\end{align}
	where $f\in\mathcal C^\infty (\RR^d, \RR^d)$ is globally Lipschitz continuous with bounded derivatives of every positive order, \(\sigma\in\RR^{d\times m}\), and \(W\) is an \(m\)-dimensional standard Wiener process. Without loss of generality, we set the terminal time $T=1$. 
	Let \(X^{N}\) denote the numerical approximation obtained by the Euler–Maruyama method. We establish the LIL for the terminal error along the geometric refinement sequence $(N_i)_{i\in\N^+}=(\alpha^i)_{i\in\N^+}$ with a positive integer $\alpha\geq 2 $ as follows 
	\begin{align*}
		\limsup_{i\to\infty}
		\frac{N_i| X(1)-X^{N_i}(1)|}
		{\sqrt{\log\log N_i}}
		=\sqrt{2\lambda_{\max} (\mathcal V_\sigma)},
		\quad\text{a.s.},
	\end{align*}
	where $ \mathcal{V}_{\sigma}:=\frac{1}{12}  \int^1_0\mathcal J(1,r) Df(X(r))\sigma \sigma^\top (Df(X(r)))^\top (\mathcal J(1,r))^\top  dr$ and $\mathcal J(t,s), t\geq s$ denotes the variational propagator.
We also establish two-sided LILs for the projections of the terminal error onto arbitrary fixed directions, identifying their sharp upper and lower asymptotic envelopes. These results determine the path-dependent asymptotic error envelope and identify \(N_i^{-1}\sqrt{\log\log N_i}\) as the exact almost sure fluctuation scale, thereby identifying the precise logarithmic correction to the first-order error scale.

The corresponding LILs also hold for the \(\theta\)-Euler–Maruyama method for every \(\theta\in[0,1]\), with the same path-dependent asymptotic constants, and extend to arbitrary fixed terminal times. These results also reveal a connection with the asymptotic error distribution, studied in, e.g., \cite{JacodProtter1998AOP, KurtzProtter1991SA, ProtterSPA2020}. Specifically, the conditional covariance matrix of the stable limit of the rescaled error, given the driving Wiener process, coincides with \(\mathcal V_\sigma\). This connection provides insight into why the constants are independent of \(\theta\), although the conditional mean of the limiting error generally depends on \(\theta\).


The proof of the LILs consists of two main steps: identifying the leading term $(Y_{N_i})_{i\in\N^+}$ in the rescaled error $N_i(X(1)-X^{N_i}(1))$ and establishing its almost sure asymptotics. 
The first step requires an error expansion that captures the oscillatory structure governing the LIL scale. Through linearization, we identify an oscillatory Skorohod integral \(Y_N\) as the leading term. 
Malliavin regularity estimates and Meyer’s inequality allow us to control the non-adapted remainder terms and show that they are negligible under the LIL normalization.
The second step therefore focuses on the almost sure asymptotics of \((Y_{N_i})_{i\in\mathbb N^+}\). The main difficulties arise from the anticipating terminal variational propagator, the generally non-Gaussian distributions of these integrals, and their dependence across refinement levels induced by the common driving Wiener process.
Isolating the terminal random factor makes martingale estimates available for the upper bound, but does not establish that fluctuations attain the predicted scale infinitely often. To obtain the matching lower bound, the key construction is a conditionally Gaussian sequence on a coarse grid whose covariance converges to the target random covariance as the grid is refined. Geometric refinement yields exponential decay of its conditional correlations, enabling a conditional Kochen–Stone argument to establish the lower bound for this sequence. We then recover the sharp bound for the original sequence by combining convergence of the coarse-grid covariance with approximation estimates that become negligible at the LIL scale as the coarse grid is refined.


	
We further compare numerical methods through their LIL constants and illustrate the theoretical results numerically.
A linear scalar example shows that the \(\theta\)-Euler–Maruyama methods can exhibit substantial differences at finite refinement levels despite sharing the same LIL constants. For a linear underdamped Langevin equation, we derive the LIL constants of the BAOAB method and show that they differ from those of the \(\theta\)-Euler–Maruyama methods despite their common strong convergence order. This analysis also illustrates how our strategy can be adapted to other numerical methods.

	The rest of this paper is organized as follows. Our main results on the almost sure asymptotics of the error process are presented in \Cref{sec_Main}, together with an outline of the proof. Section \ref{sec_preEstimates} identifies the leading term $(Y_{N_i})_{i\in\N^+}$ and controls the remainders, while Section \ref{sec_proofofThm} establishes its almost sure asymptotics.
Section~\ref{sec_Comparison} compares LIL constants and finite-refinement errors through linear examples and numerical experiments.

	\section{Almost sure asymptotics of the error process}\label{sec_Main}
	In this section, we present our main results on the almost sure asymptotics of the error processes for the Euler--Maruyama method applied to $d$-dimensional SDEs with additive noise. These results take the form of the LIL, characterizing the sharp upper envelope of the norms of the normalized errors, as well as the sharp upper and lower envelopes of their projections onto any fixed direction.

	\subsection{Notation and numerical approximation}
	We first introduce the notation used in this paper. Each Euclidean space  $\RR^d$ is equipped  with the Euclidean norm $|\cdot|$.
    For \(k\in\mathbb N^+\), positive integers \(d_0,\ldots,d_k\), a tensor \(M\in\mathbb R^{d_0\times\cdots\times d_k}\), and vectors \(u_j\in\mathbb R^{d_j}\), \(j=1,\ldots,k\), we define \(M(u_1,\ldots,u_k)\in\mathbb R^{d_0}\) by
	\[
	\bigl(M(u_1,\ldots,u_k)\bigr)_{i_0}
	:=
	\sum_{i_1=1}^{d_1}\cdots\sum_{i_k=1}^{d_k}
	M_{i_0i_1\cdots i_k}
	(u_1)_{i_1}\cdots(u_k)_{i_k},
	\qquad i_0=1,\ldots,d_0,
	\]
    where $M_{i_0i_1\cdots i_k} $ denotes the $(i_0,i_1,\ldots, i_k)$-th component of $M$.
	We denote its operator norm by
	$$
	\|M\|_{\op}
	:=\sup_{|u_1|\leq1,\ldots,|u_k|\leq1}
	|M(u_1,\ldots,u_k)|.
	$$
	The zero vector in $\RR^d$ is defined by $\mathbf 0_d$.
	For $f\in\mathcal C^k(\RR^d,\RR^d)$, we denote its $j$-th
	derivative by $D^j f$, $j=1,\ldots,k$.
	Throughout this paper, $\mathcal C_b^\infty$ denotes the space 
	of smooth functions whose derivatives of every positive order
	are bounded, while the functions themselves need not be bounded.
	We denote the It\^o integral by $\int\cdot\,dW$ and the Skorohod
	integral by $\int\cdot\,\delta W$.
	For a Malliavin differentiable matrix-valued random variable $F$,
	we write $(\mathcal D_\theta F[u])v$ for the action on $v$ of the
	Malliavin derivative of $F$ at time $\theta$ in the noise direction $u$.
	For $r\in\RR^+$, we write $[r]$ and
	$\{r\}:=r-[r]$ for its integer and fractional parts, respectively.
	Let $C$ denote a generic positive constant whose value may vary from one place to another.
	
	Next, we introduce the Euler--Maruyama method for \eqref{equ_nonliearGeneral} and the corresponding error process.
	To approximate \eqref{equ_nonliearGeneral} on $[0,1]$, we introduce
	the uniform time grid $t_k:=kh$, $k=0,\ldots,N$, with $h:=1/N$ and $N\in\N^+$.
	We also define the left and right grid-point maps by
	\[
	\eta_N(s):=t_k \quad \text{for  } t_k\leq s<t_{k+1},
	\qquad
	\zeta_N(s):=t_{k+1} \quad \text{for  } t_k<s\leq t_{k+1},
	\]
	for $k=0,\ldots,N-1$, with $\eta_N(1):=1$ and $\zeta_N(0):=0$.
	The Euler--Maruyama approximation for \eqref{equ_nonliearGeneral} is then given by
	\begin{align*}
		X_{k+1} = X_k + f(X_k)h + \sigma \Delta W_k,  \quad k=0,\ldots,N-1
	\end{align*}
	with $X_0=x_0$ and 
	$\Delta W_k:=W(t_{k+1})-W(t_k)$. 
	 Let $X^N$ 
    denote the
	continuous-time extension of the numerical solution, 
    defined by \begin{align}\label{equ_EMNon}
		X^N(t)=X_0+\int_0^t f(X^N(\eta_N(s)))ds+\sigma W(t), \quad t\in[0,1],
	\end{align}
	which satisfies $X^N(t_k)=X_k$ for $  k=0,\ldots,N$.
	Moreover, we introduce the error process
	\begin{align*}
		\mathcal{E}_N(t):=X(t)-X^N(t) = \int^t_0 \big[f(X(s))- f(X^N(\eta_N (s)))\big] ds, \quad t\in[0,1].
	\end{align*}
	
	\subsection{Main results on the almost sure asymptotics}
	In the following theorem, 
	we establish LILs along geometric mesh refinements, identifying the sharp upper and lower asymptotic envelopes of the normalized errors along any fixed direction, as well as the sharp upper asymptotic envelope of their norm.
	\begin{theorem}\label{thm_LILadditive}
		Let $(N_i)_{i\in\mathbb N^+}:=(\alpha^i)_{i\in\mathbb N^+}$ for a fixed integer \(\alpha\ge2\). Then the following results hold.
		
        \noindent(\romannumeral 1) For any given $v\in\RR^d$, 
		\begin{align}
			\limsup_{i\to \infty} \frac{N_i v^\top (X(1)-X^{N_i}(1))}{\sqrt{\log\log N_i} }&=\sqrt{2v^\top \mathcal{V}_{\sigma} v},\quad \text{a.s.}, \label{equ_mainLILsup} \\
			\liminf_{i\to \infty} \frac{N_iv^\top(X(1)-X^{N_i}(1))}{\sqrt{\log\log N_i} }&=-\sqrt{2v^\top \mathcal{V}_{\sigma}v },\quad \text{a.s.}\label{equ_mainLILinf}
		\end{align}
		where $ \mathcal{V}_{\sigma}:=\frac{1}{12}  \int^1_0\mathcal J(1,r) Df(X(r))\sigma \sigma^\top (Df(X(r)))^\top (\mathcal J(1,r))^\top  dr$ and $\mathcal J(t,s)=U(t)(U(s))^{-1}$, $t\geq s$ with $U$ solving
			\begin{align}\label{equ_UEqu}
			dU (t) = Df(X(t)) U(t)dt,\quad U(0)=I_d.
		\end{align}
		
		\noindent(\romannumeral 2) In particular, for the norm of $N_i(X(1)-X^{N_i}(1))$, 
		\begin{align}\label{equ_normLIL}
			\limsup_{i\to \infty} \frac{|N_i(X(1)-X^{N_i}(1))|}{\sqrt{\log\log N_i} }=\sqrt{2\lambda_{\max}(\mathcal V_\sigma)},\quad \text{a.s.},
		\end{align}
		where $\lambda_{\max}(\mathcal V_{\sigma})$ represents the largest eigenvalue of $\mathcal V_{\sigma}$. 
	\end{theorem}
 
 The almost-sure error asymptotics established in \Cref{thm_LILadditive} remain valid for the \(\theta\)-Euler–Maruyama method with the same path-dependent asymptotic constants and extend to arbitrary fixed terminal times. 
	\begin{rem} \label{thm_thetamethod}

\noindent(\romannumeral 1)	For any \(\theta\in[0,1]\), consider the \(\theta\)-Euler–Maruyama method applied to \eqref{equ_nonliearGeneral}:
\begin{align} \label{equ_thetaEM}
X_{k+1}^{\theta}
=X_k^\theta
+h\big[(1-\theta)f(X_k^\theta)
+\theta f(X_{k+1}^\theta)\big]+\sigma\Delta W_k,
\qquad k=0,\ldots,N-1.
\end{align}
Denote by \(X^{N,\theta}(1)=X_N^\theta\) the numerical approximation at the terminal time with step size \(h=1/N\). Then, for every \(\theta\in[0,1]\), the LILs in Theorem~\ref{thm_LILadditive} remain valid for the terminal error \(\mathcal E_{N,\theta}(1):=X(1)-X^{N,\theta}(1)\), with the same path-dependent asymptotic constants as those of the Euler–Maruyama method. The proof follows the same approach as that of Theorem~\ref{thm_LILadditive} and is deferred to Section~\ref{subsec_prooftheta}.
	
\noindent(\romannumeral 2)	Keep $h_i=N_i^{-1}$.	For any fixed terminal time \(T>0\), the corresponding LILs hold for the terminal error \(X(T)-X^{N_i}(T)\), with \(\mathcal V_\sigma\) replaced by
        \[
        \mathcal{V}_{T,\sigma}:=\frac{1}{12}  \int^T_0\mathcal J(T,r) Df(X(r))\sigma\sigma^\top  (Df(X(r)))^\top (\mathcal J(T,r))^\top  dr.
        \]
	\end{rem}

The path-dependent LIL constants also admit a natural interpretation through the asymptotic error distribution. In particular, their connection with the conditional covariance of the limiting error provides insight into why these constants are independent of \(\theta\).
    
	\begin{rem} \label{rem_asymptoticCov}
    
        
        For the \(\theta\)-Euler–Maruyama method applied to \eqref{equ_nonliearGeneral}, 
        the normalized error process
		$\mathcal U_{N,\theta}:=N\mathcal E_{N,\theta}$  stably converges in distribution in
		$\mathcal C([0,1];\RR^d)$ to the limiting process $\mathcal U_\theta$, which satisfies
		\begin{align}\label{equ_AED}
			d\,\mathcal U_\theta(t)
			&= Df(X(t))\,\mathcal U_\theta(t)dt
			+ \big(\frac{1}{2}-\theta\big) \Big(
			Df(X(t))f(X(t))
			+ \frac{1}{2}\sum_{j=1}^{m}
			D^2f(X(t))(\sigma e_j,\sigma e_j)
			\Big)\,dt \nn\\
			&\quad
			+ \big(\frac{1}{2}-\theta\big)Df(X(t))\sigma dW(t)
			+ \frac{1}{\sqrt{12}}Df(X(t))\sigma dB(t),
			\qquad  \mathcal U_\theta(0)=0,
		\end{align}
		where $B$ is an $m$-dimensional standard Wiener process
		independent of $W$, and ${(e_j)}_{j=1}^m$ is the standard
		basis of $\RR^m$; see e.g. \cite{KurtzProtter1991SA} for the study of asymptotic error distributions.
		Since \eqref{equ_AED} is linear in $\mathcal U_\theta$, we obtain
		the following explicit representation
		\begin{align*}
				\mathcal U_\theta(t) &=\big(  \frac{1}{2}-\theta\big) U(t)\int^t_0 U(s)^{-1}\Big(Df(X(s))f(X(s)) +\frac{1}{2}\sum^{m}_{j=1} D^2f(X(s)) (\sigma e_j,\sigma e_j) \Big)ds \\
			&\quad+ \big(  \frac{1}{2}-\theta\big)  U(t)\int^t_0 U(s)^{-1} Df(X(s))\sigma dW(s)+ \frac{1}{\sqrt{12}} U(t)\int^t_0 U(s)^{-1} Df(X(s)) \sigma dB(s).
		\end{align*}
	Therefore,  the asymptotic conditional
	covariance matrix coincides with $\mathcal V_{\sigma}$, i.e., 
    $$
			\operatorname{Cov}(\,	\mathcal U_\theta(1)\,|\, \mathcal F_W ) = \mathcal{V}_{\sigma},\quad \mathcal F_W:=\sigma\{ W(t): 0\leq t\leq 1\}
	$$
	This reveals the connection between the asymptotic error distribution
	and the almost sure asymptotics. The conditional mean  $\E[\,	\mathcal U_\theta(1)\, |\, \mathcal F_W] $ becomes
	negligible under the LIL normalization, while the conditional
	covariance arising from the independent Wiener process $B$
	determines the sharp envelopes of the almost sure error fluctuations.
	This also provides an intuitive explanation of why the LILs for the $\theta$-Euler--Maruyama methods have the same asymptotic constants for all $\theta\in[0,1]$. 
	\end{rem}

    \subsection{Proof of \Cref{thm_LILadditive}} 
	We now turn to the proof of our main theorem. 
     For any strictly increasing sequence of positive integers
$(N_i)_{i\in\N^+}$, we identify the leading term in the rescaled error $N_i\mathcal E_{N_i}(1)$ as
	\begin{align}\label{equ_YNexpress}
		Y_{N_i}:=  \int^1_0  \mathcal J(1,r) Df(X(r))\big(\frac 12-\{N_i r\} \big) \sigma \delta W(r).
	\end{align}
    The following proposition shows that the remainder
$N_i\mathcal E_{N_i}(1)-Y_{N_i}$ vanishes almost surely
under the LIL normalization.
The proof is deferred to \Cref{sec_preEstimates}.
	\begin{prop}\label{prop_remainder}
		Let ${(N_i)}_{i\in\N^+}$ be a strictly increasing sequence of positive integers. Then we have 
		\begin{align*}
			\limsup_{i\to \infty} \frac{ |N_i {\mathcal E}_{N_i} (1)-Y_{N_i} |}{\sqrt{\log\log N_i}} = 0,\quad \text{a.s.}
		\end{align*} 
	\end{prop}
	
    Having established that the remainder is negligible, we now focus on the leading term $Y_{N_i}$.
 	The following proposition characterizes the sharp upper and lower
 	envelopes of its  normalized projections onto any fixed
 	direction. 
 	The proof of the proposition is deferred to \Cref{sec_proofofThm}.
 	
		\begin{prop}\label{prop_LILofY}
		Let $v\in\RR^d$ and $(N_i)_{i\in\mathbb N^+}:=(\alpha^i)_{i\in\mathbb N^+}$ for a fixed integer \(\alpha\ge2\). Then it holds that
		\begin{align*}
			\limsup_{i\to \infty} \frac{ v^\top Y_{N_i}}{\sqrt{\log\log N_i} }=\sqrt{2v^\top \mathcal{V}_{\sigma} v},\quad \text{a.s.}\quad \text{and} \quad \liminf_{i\to \infty} \frac{ v^\top Y_{N_i}}{\sqrt{\log\log N_i} }=-\sqrt{2v^\top \mathcal{V}_{\sigma} v},\quad \text{a.s.},
		\end{align*}
		where $\mathcal{V}_{\sigma}$ is introduced in \Cref{thm_LILadditive}. 
	\end{prop}
    With these results, we now present the proof of Theorem~\ref{thm_LILadditive}.
	\begin{proof}[Proof of Theorem~\ref{thm_LILadditive}]
    \noindent (\romannumeral 1)
	Combining Propositions~\ref{prop_remainder} and \ref{prop_LILofY}, for any $v\in\RR^d$, we obtain
	\begin{align*}
				\limsup_{i\to \infty} \frac{N_i v^\top (X(1)-X^{N_i}(1))}{\sqrt{\log\log N_i} } &\leq 		\limsup_{i\to \infty} \frac{v^\top Y_{N_i} }{\sqrt{\log\log N_i} } + 	\limsup_{i\to \infty} \frac{ |v||N_i {\mathcal E}_{N_i} (1)-Y_{N_i} |}{\sqrt{\log\log N_i}} \\
				&\leq \sqrt{2v^\top \mathcal{V}_{\sigma} v},\quad \text{a.s.}
	\end{align*}
	and 
		\begin{align*}
		\limsup_{i\to \infty} \frac{N_i v^\top (X(1)-X^{N_i}(1))}{\sqrt{\log\log N_i} } &\geq 		\limsup_{i\to \infty} \frac{v^\top Y_{N_i} }{\sqrt{\log\log N_i} } -	\limsup_{i\to \infty} \frac{ |v||N_i {\mathcal E}_{N_i}(1) -Y_{N_i} |}{\sqrt{\log\log N_i}} \\
		&\geq \sqrt{2v^\top \mathcal{V}_{\sigma} v},\quad \text{a.s.}
	\end{align*}
	These two inequalities complete the proof of the $\limsup$ identity in  \Cref{thm_LILadditive} (\romannumeral 1), while the $\liminf$ identity follows by the same argument. 
	
	 \noindent (\romannumeral 2) We first claim that by the compactness of the unit sphere $\mathbb S^{d-1}$ in $\RR^d$, for every 
	$\varepsilon\in(0,1)$, there exists a finite $\varepsilon$-net $\{v_1,\ldots, v_m \}$ of $\mathbb S^{d-1}$ such that
	\begin{align}\label{equ_vectortonorm}
		(1-\varepsilon)|z| \leq \max_{1\leq j\leq m} v_j^\top z,\quad z\in\RR^d.
	\end{align}
In fact, for  $z=\mathbf 0_d$, \eqref{equ_vectortonorm} obviously holds.  For $z\in\RR^d\backslash \{\mathbf 0_d\}$, there exists $j\in\{1,\ldots,m\}$ such that $|\frac{z}{|z|}-v_j|\leq \varepsilon$, which leads to 
	\begin{align*}
		|z| = \Big(\frac{z}{|z|}\Big)^\top z \leq v_j^\top z + \Big| \frac{z}{|z|}-v_j\Big| |z|   \leq v_j^\top z + \varepsilon |z|
	\end{align*} and thereby proves \eqref{equ_vectortonorm}.

    Next we prove \eqref{equ_normLIL} based on the directional LIL in (\romannumeral 1) and \eqref{equ_vectortonorm}. 
    For each integer \(n\ge2\), let \(\mathscr D_n\) be a finite \(1/n\)-net of \(\mathbb S^{d-1}\) satisfying \(\eqref{equ_vectortonorm}\) with \(\varepsilon=1/n\). Set
$
\mathscr D=\cup_{n=2}^\infty\mathscr D_n.
$
Then \(\mathscr D\) is a countable dense subset of \(\mathbb S^{d-1}\).
      It follows from \Cref{thm_LILadditive} (\romannumeral 1) that there exists $\Omega_0\subset \Omega$ with $\PP( \Omega_0)=1$ such that \eqref{equ_mainLILsup}--\eqref{equ_mainLILinf} hold for all $v\in\mathscr D$.   
	For any $\omega\in\Omega_0$, it follows from \eqref{equ_vectortonorm} with the net $\mathscr D_n$ that
	\begin{align*}
		\limsup_{i\to \infty} \frac{ |N_i  (X(1,\omega)-X^{N_i}(1,\omega))|}{\sqrt{\log\log N_i} } \leq \frac{1}{1-1/n} \max_j \sqrt{2v_j^\top \mathcal{V}_{\sigma}(\omega) v_j} \leq \frac{1}{1-1/n} \sqrt{2\lambda_{\max}(\mathcal{V}_{\sigma}(\omega))}.
	\end{align*}
	Letting $n\to \infty$ gives the upper bound of \eqref{equ_normLIL}. 
	The lower bound follows from the fact that for any $\omega\in \Omega_0$,
	\begin{align*}
		\limsup_{i\to \infty} \frac{ |N_i  (X(1,\omega)-X^{N_i}(1,\omega))|}{\sqrt{\log\log N_i} } &\geq \sup_{v\in\mathscr D}\limsup_{i\to \infty} \frac{ N_i v^\top (X(1,\omega)-X^{N_i}(1,\omega))}{\sqrt{\log\log N_i} }  \\
        &=\sup_{v\in\mathscr D} \sqrt{2v^\top \mathcal V_{\sigma} (\omega)v} 
		 =\sqrt{2\lambda_{\max}(\mathcal{V}_{\sigma}(\omega))},
	\end{align*}
	which ends the proof. 
	\end{proof}
    \section{Identification of the leading term $Y_{N_i}$} \label{sec_preEstimates}
	In this section, we identify the leading term $Y_{N_i}$ and prove \Cref{prop_remainder}. Specifically, we first derive a uniform rough upper bound for $N_i\mathcal{E}_{N_i}$, and then drop the higher-order terms in the Taylor expansion by showing the equivalence of the linearization process. Finally, we identify $Y_{N_i}$ from the linearization, proving that all the remainders vanish almost surely after division by \(\sqrt{\log\log N_i}\).

	\subsection{Uniform rough upper bound of the error process $\mathcal{E}_{N_i}$}
	We provide a rough upper bound for $\mathcal{E}_{N_i}$ in the following proposition, which shows that $N_i\mathcal{E}_{N_i}$ vanishes under any positive-order normalization of $N_i$.
	\begin{prop} \label{prop_roughupperbound}
		For any given $\varepsilon>0$ and strictly increasing
		sequence of positive integers $\{N_i\}_{i\in\N^+}$, it holds that
		\begin{align*} 
			\lim_{i\to \infty} \frac{N_i	\sup_{t\in[0,1]} | 	\mathcal{E}_{N_i}(t) | }{(N_i)^{\varepsilon}} =0, \quad \text{a.s.}
		\end{align*} 
	\end{prop}
	\begin{proof}
			It is well-known that for any given $p\geq 2$, there exists a constant $C_p>0$ independent of $N$ such that
		\begin{align}\label{equ_UniformLpE}
			\Big( \E\Big[ \sup_{t\in[0,1]}\big|\mathcal E_N(t)\big|^p \Big]\Big)^{\frac 1p} 	\leq
			\frac{C_p}{N}.
		\end{align}

	Let $\varepsilon>0$ and choose $p\geq2$ such that $p\varepsilon>1 $. For any given strictly increasing sequence of positive integers $\{N_i\}_{i\in\N^+}$ and any given 
		$\tilde{\varepsilon}>0$, the Chebyshev inequality gives
		\begin{align*}
			\mathbb P\Big\{
			\frac{	N_i \sup_{t\in[0,1]}|\mathcal E_{N_i}(t)|	}{	(N_i)^\varepsilon	}	>\tilde{\varepsilon}
			\Big\}&\leq
			\frac{	\mathbb E\big[
				\sup_{t\in[0,1]}|\mathcal E_{N_i}(t)|^p
				\big]		}{	(\tilde{\varepsilon})^p(N_i)^{p(\varepsilon-1)}	}
			\leq
			\frac{C_{p,\tilde{\varepsilon}}}
			{(N_i)^{p\varepsilon}}.
		\end{align*}
		Note that $
		\sum_{i=1}^\infty
		\frac{1}
		{(N_i)^{p\varepsilon}}
		<\infty
		$ holds for $N_i\geq i$ and $p\varepsilon>1$.
		Combining this with  the  Borel--Cantelli lemma, we derive 
		\begin{align*}
			\PP \Big\{
			\frac{	N_i \sup_{t\in[0,1]}|\mathcal E_{N_i}(t)|	}{	(N_i)^\varepsilon	}	> \tilde{\varepsilon} \quad \text{i.o.} 
			\Big\}=0
		\end{align*}
		for any $\tilde{\varepsilon}>0$. Taking  $\tilde{\varepsilon} \in \mathbb {Q}^+$ such that $\tilde{\varepsilon} \to 0$ therefore yields 
		$$
		\lim_{i\to\infty}
		\frac{
			N_i \sup_{t\in[0,1]}|\mathcal E_{N_i}(t)|
		}{
			(N_i)^\varepsilon
		}	=0,
		\quad\text{a.s.}
		$$
		for any $\varepsilon>0$. The proof is complete.
	\end{proof}
	\subsection{Linearization of the error process} We first derive a linearized approximation of the error process \(\mathcal E_N\) using the Taylor expansion and then establish their asymptotic equivalence under the LIL normalization.
	\begin{align}\label{equ_ENtaylor}
		\mathcal{E}_N(t) &= \int_0^t [f(X(s)) - f(X^N(\eta_N(s)))] ds \nn\\
		&= \int_0^t\Big[  Df(X(s))\big(X(s) - X^N(\eta_N(s))\big) - \bar{\xi}_N(s) \big(X(s) - X^N(\eta_N(s)),X(s) - X^N(\eta_N(s)) \big) \Big]ds\nn \\
		&=  \int_0^t  Df(X(s))\big(\mathcal{E}_N(s)+ f(X^N(\eta_N(s)))(s - \eta_N(s)) + \sigma(W(s) - W(\eta_N(s)))\big) ds \nn \\
		&\quad   -\int_0^t  \bar{\xi}_N(s) \Big(\mathcal{E}_N(s)+ f(X^N(\eta_N(s)))(s - \eta_N(s)) + \sigma(W(s) - W(\eta_N(s))), \nn \\ 
		&\qquad \qquad \mathcal{E}_N(s)+ f(X^N(\eta_N(s)))(s - \eta_N(s)) + \sigma(W(s) - W(\eta_N(s)))\Big) ds,
	\end{align}
	where we use
	\begin{align*}
		X(s) - X^N(\eta_N(s))&= 	\mathcal{E}_N(s)+ 
		X^N(s) - X^N(\eta_N(s)) \\
		&=	\mathcal{E}_N(s)+  f(X^N(\eta_N(s)))(s - \eta_N(s)) + \sigma(W(s) - W(\eta_N(s)))
	\end{align*}
	and
	$\bar{\xi}_N(s) (u,v) :=\int_0^1 \theta 
	D^2f\big( \theta
	X(s)+ (1-\theta) X^N(\eta_N(s))
	\big) (u,v) d\theta$ for $u,v\in\RR^d$. 
	According to \eqref{equ_ENtaylor}, 
	we construct an auxiliary process $\widetilde{\mathscr{E}}_N$ satisfying the following linear equation
	\begin{align} \label{equ_VN}
		d \widetilde{\mathscr{E}}_N(t) &=Df(X(t))\widetilde{\mathscr{E}}_N(t) dt +Df(X(t)) \big[ f(X^N(\eta_N(t)))(t - \eta_N(t)) + \sigma(W(t) - W(\eta_N(t)))\big] dt
	\end{align}
    with $\widetilde{\mathscr E}_N(0)=0$.
	Recall that $\mathcal J(t,s)=U(t) (U(s))^{-1}, \, t\geq s\ge  0$ with $U(t)$ solving 
	\eqref{equ_UEqu}. Then $\widetilde{\mathscr{E}}_N$ admits the following explicit expression
	\begin{align} \label{equ_exactVn}
		\widetilde{\mathscr{E}}_N(t)  
		&= \int^t_0 \mathcal J(t,s)  Df(X(s)) \big[f(X^N(\eta_N(s)))(s - \eta_N(s)) + \sigma(W(s) - W(\eta_N(s)))\big] ds\nn \\
		&=\int^t_0  \mathcal J(t,s)  Df(X(s)) f(X^N(\eta_N(s)))(s - \eta_N(s))  ds \nn \\
		&\quad +\int^t_0  \mathcal J(t,s) Df(X(s))\sigma  (W(s) - W(\eta_N(s))) ds.
	\end{align}
	
	By taking the derivative on both sides of the identity $(U(t))^{-1} U(t) =I_d$ with respect to $t$ and using \eqref{equ_UEqu}, we have
	\begin{align*}
		\big(d (U(t))^{-1}\big) U(t) +(U(t))^{-1} dU(t)  =  \big(d (U(t))^{-1}\big) U(t) +(U(t))^{-1} Df(X(t))U(t) dt =0,
	\end{align*}
	which yields
	\begin{align}\label{equ_UinverseEqu}
		d (U(t))^{-1} = -(U(t))^{-1} Df(X(t)) dt, \quad (U(0))^{-1}=I_d.
	\end{align}
	
	Before showing the equivalence between the original process $\mathcal E$ and the
	linearized error process $\widetilde{\mathscr E}$, we introduce some  properties of $\mathcal J$, $U$, and $U^{-1}$, which will be frequently used in the following content.
	
	\textbf{Boundedness of  $\| U(t)\|_{\op}$,  $\| (U(t))^{-1}\|_{\op}$, and $\|\mathcal J(t,s)\|_{\op}$ for $0\leq s\leq t\leq 1$.} For any $v\in\RR^d$ and $s\in[0,1]$,  it follows from \eqref{equ_UEqu} and  the uniform boundedness $\|Df(\cdot)\|_{\op}\leq C$ that
	\begin{align*}
		d| U(t)v|^2 = 2(U(t)v)^\top Df(X(t)) U(t)v dt \leq C |U(t)v|^2 dt. 
	\end{align*}
	Applying the Gr\"onwall inequality leads to $|U(t)v|\leq C|v|$, which further yields $\| U(t)\|_{\op}\leq C$ for $t\in[0,1]$ with $C$ independent of $t$ and $s$. Following a similar procedure, we also have  $\| (U(t))^{-1}\|_{\op}\leq C$ for $t\in[0,1]$  and then $\| \mathcal J(t,s)\|_{\op}\leq \| U(t)\|_{\op}\| (U(s))^{-1}\|_{\op}\leq C$ for $t\in[s,1]$ with $C$ independent of $s,t$. Therefore, we have 
	\begin{align*}
		\sup_{t\in[0,1]} \| U(t)\|_{\op}\leq C, \quad  \sup_{t\in[0,1]}\| (U(t))^{-1}\|_{\op}  \leq C \quad \text{and}\quad  \sup_{0\leq s\leq t\leq 1}\|\mathcal J(t,s)\|_{\op} \leq C.
	\end{align*}
	
	\textbf{First-order temporal continuity of  $U$, $U^{-1}$, and $\mathcal J$.}
	Recall that $\mathcal J(1,r)=U(1)(U(r))^{-1}$ with $U$ satisfying \eqref{equ_UEqu}. 
	It follows from \eqref{equ_UinverseEqu} and the uniform boundedness $\sup_{r\in[0,1]}\|(U(r))^{-1}\|_{\op}\leq C$, $\sup_{x\in\RR^d}\|Df(x)\|_{\op}\leq C$ that 
	\begin{align} \label{equ_JTemporalHolder}
		\| \mathcal J(t,s) - \mathcal J(t,r)\|_{\op} &=  \| U(t)( (U(s))^{-1}-(U(r))^{-1}) \|_{\op}  \nn \\
		&\leq C \|  (U(s))^{-1}-(U(r))^{-1}\|_{\op} \nn \\
		&\leq C \Big\| \int^r_s(U(z))^{-1} Df(X(z)) dz \Big\|_{\op} \leq C|r-s|, \quad 0\leq r,s\leq t\leq 1.
	\end{align}
	Similarly, we also have $\|\mathcal J(r,s) - \mathcal J(t,s)\|_{\op}\leq C|r-t|$ for $0\leq s\leq r, t\leq 1$, $\|  U(s)-U(r)\|_{\op}\leq C|r-s|$ and $\|  (U(s))^{-1}-(U(r))^{-1}\|_{\op}\leq C|r-s|$ for $0\leq r,s\leq 1$.

	We finally emphasize that $U(t)$ cannot generally be represented
	as the exponential of $\int_0^t Df(X(r))\,dr$ for $d>1$, since the matrices
	$Df(X(t))$ at different times need not commute.
	However, explicit representations are available when $d=1$
	or $Df(X(t))$ is a constant matrix. Specifically, for $d=1$,
	\begin{align*}
		U(t)=\exp\big\{\int_0^t Df(X(r))dr\big\} \quad \text{and}\quad \mathcal J (t,s)= \exp\big\{  \int^t_s Df(X(r)) dr \big\}.
	\end{align*}
	 If $Df(X(t))\equiv A$ for some constant matrix $A$, then
	 	\begin{align*}
	 U(t)=e^{At} \quad \text{and}\quad \mathcal J (t,s)= e^{A(t-s)}.
	 \end{align*}

    These estimates allow us to establish the asymptotic equivalence between the original and linearized error processes under the LIL normalization.
	
	\begin{prop}\label{prop_linearization}
		For any
		 strictly increasing	sequence of positive integers $(N_i)_{i\in\N^+}$,
        it holds that 
		\begin{align*}
			\lim_{i\to \infty} \frac{ N_i  \sup_{t\in[0,1]}	| 	\mathcal{E}_{N_i}(t)-	\widetilde{\mathscr{E}}_{N_i}(t)  |} {\sqrt{\log \log N_i}}=0, \quad \text{a.s.} 
		\end{align*}
	\end{prop} 
	\begin{proof}
		Combining \eqref{equ_ENtaylor} and \eqref{equ_VN}, we have 
		\begin{align*}
			\mathcal{E}_N(t)-	\widetilde{\mathscr{E}}_N(t)  
			&=  \int_0^t  Df(X(s))\big( \mathcal{E}_N(s) - 	\widetilde{\mathscr{E}}_N(s) \big) ds \\
			&\quad -\int_0^t \bar{\xi}_N(s)\big[\mathcal{E}_N(s)+ f(X^N(\eta_N(s)))(s - \eta_N(s)) + \sigma(W(s) - W(\eta_N(s))),\\
			&\qquad \qquad \mathcal{E}_N(s)+ f(X^N(\eta_N(s)))(s - \eta_N(s)) + \sigma(W(s) - W(\eta_N(s))) \big] ds.
		\end{align*}
		Using the uniform boundedness $\|Df\|_{\op}\leq C$ and $\| \bar\xi_N\|_{\op} \leq \|D^2f\|_{\op} \leq C$,
		we derive
		\begin{align*}
			| 	\mathcal{E}_N(t)-	\widetilde{\mathscr{E}}_N(t)  | &\leq  C  \int_0^t 	| 	\mathcal{E}_N(s)-	\widetilde{\mathscr{E}}_N(s)  | ds \\
			&\quad +  C \int_0^t \big|\mathcal{E}_N(s)+ f(X^N(\eta_N(s)))(s - \eta_N(s)) + \sigma(W(s) - W(\eta_N(s)))\big|^2 ds,
		\end{align*}
		upon which the  Gr\"onwall inequality gives
		\begin{align}\label{equ_Ewhole}
			&\quad 	\sup_{t\in[0,1]} | 	\mathcal{E}_N(t)-	\widetilde{\mathscr{E}}_N(t)  | \nn\\
			&\leq C    \int_0^1 \big|\mathcal{E}_N(s)+ f(X^N(\eta_N(s)))(s - \eta_N(s)) + \sigma(W(s) - W(\eta_N(s)))\big|^2 ds\nn \\
			&\leq  C   \Big( \int_0^1 |\mathcal{E}_N(s)|^2 ds+ \int^1_0  |f(X^N(\eta_N(s)))|^2(s - \eta_N(s))^2 ds+  \|\sigma\|^2_{\op}\int^1_0 |W(s) - W(\eta_N(s)) |^2 ds \Big).
		\end{align}
		
	In view of \eqref{equ_Ewhole}, it remains to show that each term on the right-hand side, after normalization, converges to zero almost surely.
		First, it follows from \Cref{prop_roughupperbound} that for any $\varepsilon\in(0,\frac 12)$,
		\begin{align*}
			\limsup_{i\to \infty} \frac{ N_i \int^1_0  |\mathcal{E}_{N_i} (s)|^2 ds}{\sqrt{\log \log N_i}}&\leq 	\limsup_{i\to \infty} \frac{ N_i (\sup_{0\leq s\leq 1}  |\mathcal{E}_{N_i} (s)|)^2 }{\sqrt{\log \log N_i}}  \\
			&\leq \limsup_{i\to \infty}  \Big(	 \frac{ N_i \sup_{0\leq s\leq 1}  |\mathcal{E}_{N_i} (s)| }{ (N_i)^{\varepsilon}} \Big)^2 \Big(\frac{ 1 }{ (N_i)^{1-2\varepsilon}} \Big)\\ 
			&\leq \Big(	\limsup_{i\to \infty} \frac{ N_i \sup_{0\leq s\leq 1}  |\mathcal{E}_{N_i} (s)| }{ (N_i)^{\varepsilon} } \Big)^2 \Big(\limsup_{i\to \infty} \frac{ 1 }{ (N_i)^{1-2\varepsilon}} \Big)=0, \quad \text{a.s.}
		\end{align*}
		Second, combining the Lipschitz property of $f$ with the identity $ \int^1_0  (s - \eta_N(s))^2 ds=\frac{1}{3N^2}$, we have
		\begin{align*}
			\E\Big[ \Big|	N_i \int^1_0  |f(X^{N_i}(\eta_{N_i}(s)))|^2(s - \eta_{N_i}(s))^2 ds\Big|^2 \Big] \leq C \frac{1}{(N_i)^2}\Big(1+ \E\big[ \sup_{s\in[0,1]} |X^{N_i}(s)|^{4} \big] \Big)\leq C \frac{1}{(N_i)^2}, 
		\end{align*}
		where the uniform moment boundedness for $X^{N} $ can be derived using standard arguments (see e.g. \cite{Milstein2004book}). Then for any $\varepsilon>0$, it holds that
		\begin{align*}
			&\quad \sum^{\infty}_{i=1}	\PP \Big\{\frac{ N_i  \int^1_0  |f(X^{N_i}(\eta_{N_i}(s)))|^2(s - \eta_{N_i}(s))^2 ds}{\sqrt{\log \log N_i}} > \varepsilon \Big\} \\
			&\leq  \sum^{\infty}_{i=1}\frac{ \E\big[ \big| N_i  \int^1_0  |f(X^{N_i}(\eta_{N_i}(s)))|^2(s - \eta_{N_i}(s))^2 ds \big|^2 \big] }{ \varepsilon^2 (\sqrt{\log\log N_i})^2   } 
			\leq C  \sum^{\infty}_{i=1}\frac{1}{ (N_i)^2  (\sqrt{\log\log N_i})^2 } <\infty, 
		\end{align*}
		which implies
		\begin{align*}
			\limsup_{i\to \infty} \frac{ N_i  \int^1_0  |f(X^{N_i}(\eta_{N_i}(s)))|^2(s - \eta_{N_i}(s))^2 ds}{\sqrt{\log \log N_i}}=0,  \quad \text{a.s.}
		\end{align*}
		Finally, we prove that $\lim_{i \to \infty} \frac{ N_i  \int^1_0|W(s) - W(\eta_{N_i}(s))|^2 ds}{\sqrt{\log \log N_i}}=0$ almost surely. 
		Since 
		\begin{align*}
			J_N:= \int_0^1 |W(s) - W(\eta_N(s))|^2 ds
		\end{align*}
		satisfies
		\begin{align*}
			\E[ NJ_N] &= N \sum_{k=0}^{N-1} \int_{\frac kN}^{\frac{k+1}{N}} \E\big[ \big|W(s) - W(\frac{k}{N})\big|^2 \big] ds =  Nm \sum_{k=0}^{N-1} \int_{\frac kN}^{\frac{k+1}{N}}  (s- \frac kN) ds =  \frac m2,
		\end{align*}
		we define the centered random variable on each small interval  by
		$$
		J_{N,k}:=N \int_{\frac kN}^{\frac{k+1}{N}}  \big|W(s) - W(\frac{k}{N})\big|^2  ds-\frac{m}{2N} 
		$$ 
		for $k=0,\ldots, N-1$.
	For each fixed \(N\), the random variables \((J_{N,k})_{k=0}^{N-1}\) are independent and centered, since they depend on Wiener increments over disjoint intervals. Hence,
    \begin{align}\label{equ_E4_N}
			\E\big[ |NJ_N -\frac m2 |^4\big] &=  \E\big[|\sum_{k=0}^{N-1}  J_{N,k}|^4 \big]  =\sum_{k=0}^{N-1}  \E \big[ | J_{N,k}|^4 \big] + 3 \sum_{k\neq  j}  \E \big[ | J_{N,k}|^2\big]  \E \big[ | J_{N,j}|^2\big]. 
		\end{align}
		For any $p\geq 2$, we have
		\begin{align*}
			\E [ | J_{N,k}|^p ]  &\leq  C_p \Big(N^p \E\Big[ \Big| \int^{\frac{k+1}{N}}_{\frac kN} \big|W(s)-W(\frac kN )\big|^{2} ds \Big|^p  \Big]  + \frac{m^p}{2^p N^p}\Big) \\
			&\leq C_p N^p  \Big( \int^{\frac{k+1}{N}}_{\frac kN} \E\big[\big|W(s)-W(\frac kN )\big|^{2p}\big]ds \Big)  \Big(\frac{1}{N} \Big)^{p-1} +C_{p}  \frac{1}{ N^p} \\
			&\leq C_{p} N \big(\int^{\frac 1N}_0 s^{p} ds\big) +C_{p} \frac{1}{ N^p} \leq  C_{p}  \frac{1}{ N^p},
		\end{align*}
		where  we use the fact that 
		$
		\E[|W(s)-W(\frac kN )|^{2p}]\leq C_{p} (s-\frac kN)^p.
		$
		Therefore, taking $p=2$ and $4$, respectively, yields
		$
		\E [ | J_{N,k}|^2 ] \leq  C\frac{1}{N^2} 
		$
		and
		$
		\E [ | J_{N,k}|^4 ] \leq C \frac{1}{N^4} 
		$
		for any $k=0,\ldots, N-1$. 
		Plugging these two estimates into \eqref{equ_E4_N} gives
		\begin{align*}
			\E\big[ |NJ_N -\frac m2 |^4\big]  \leq C\frac{1}{N^3}+ CN(N-1) \frac{1}{N^4} \leq C\frac{1}{N^2}. 
		\end{align*}
		Hence, 
		for any given $\varepsilon>0$, the Chebyshev inequality yields 
		\begin{align*}
			\sum^\infty_{N=1}	\PP\{ | NJ_N -\frac m2| >\varepsilon\} \leq 	\sum^\infty_{N=1} \frac{\E[| NJ_N -\frac m2|^4 ]}{\varepsilon^4} \leq C	\sum^\infty_{N=1} \frac{1}{\varepsilon^4N^2}<\infty.
		\end{align*}
		Letting $\varepsilon\in\mathbb Q^+$ such that $\varepsilon\to 0$, the Borel--Cantelli lemma shows that 
		\begin{align*}
			\lim_{N\to \infty} N \int_0^1 \big|W(s) - W(\eta_N(s))\big|^2 ds =\frac m2, \quad \text{a.s.},
		\end{align*}
		which immediately implies 
		\begin{align*}
			\lim_{i \to \infty} \frac{N_i \int_0^1 \big|W(s) - W(\eta_{N_i}(s))\big|^2 ds}{\sqrt{\log \log N_i}} =0, \quad \text{a.s.}
		\end{align*}
		because ${(N_i)}_{i\in \N^+}$ is a subsequence of the positive integer sequence ${(N)}_{N\in\N^+}$. 
		
		Combining the preceding estimates, we end the proof.
	\end{proof}
	\subsection{Identification of the leading term $Y_{N_i}$ from $N_i\widetilde{\mathscr{E}}_{N_i}(1)$} 
	In this subsection, we identify $Y_{N_i}$ defined in
\eqref{equ_YNexpress} as the leading term of
$N_i\widetilde{\mathscr E}_{N_i}(1)$ under the LIL normalization.
    To this end, we first derive a decomposition that isolates
$Y_{N_i}$ and then show that the remaining terms vanish as $i\to \infty$ after division by $\sqrt{\log\log N_i}$.

    We begin by decomposing the stochastic term in
\eqref{equ_exactVn}.  Define
		\begin{align}\label{eq:defg}
		g(s):=   \mathcal J(1,s) Df(X(s)), \quad s\in[0,1].
        \end{align}
        Applying the product formula for the divergence operator (see e.g. \cite[Proposition 3.2.4]{Nualart2018book}) and the Fubini theorem yields
		\begin{align}\label{equ_FubiniMain}
			&\quad	\int^1_0  \mathcal J(1,s)  Df(X(s)) \sigma (W(s) - W(\eta_N(s))) ds \nn \\
			&=  \int^1_0 g(s) \Big( \int^s_{\eta_N(s)} \sigma dW(r)\Big) ds \nn \\
			&= \int^1_0  \int^s_{\eta_N(s)}  g(s)  \sigma \delta W(r) ds +  \sum^m_{j=1} \int^1_0  \int^s_{\eta_N(s)}  ( \mathcal D_rg(s) [e_j] ) (\sigma e_j) drds   \nn \\
			&= \int^1_0  \int^{\zeta_N(r)}_{r} g(s) \sigma ds \delta W(r)+  \sum^m_{j=1} \int^1_0  \int^s_{\eta_N(s)}  ( \mathcal D_rg(s) [e_j] ) (\sigma e_j) drds  \nn \\
			&=  \int^1_0 g(r) (\zeta_N(r)-r)  \sigma \delta W(r) + 
		 \int^1_0  \int^{\zeta_N(r)}_r \big(g(s)-g(r)\big)  \sigma ds\delta W(r)\nn \\
		 &\quad + \sum^m_{j=1} \int^1_0  \int^s_{\eta_N(s)}  ( \mathcal D_rg(s) [e_j] ) (\sigma e_j) drds \nn \\
		&	=: \uppercase\expandafter{\romannumeral 1}_N+\uppercase\expandafter{\romannumeral 2}_N+\uppercase\expandafter{\romannumeral 3}_N ,
		\end{align}
		where $e_j\in\RR^m $ is the $j$-th standard basis vector. 	
        The leading term $Y_{N_i}$ arises from $N_i\uppercase\expandafter{\romannumeral 1}_{N_i}$, as shown by the following decomposition:
		\begin{align} \label{equ_Mainpart}
			N_i \uppercase\expandafter{\romannumeral 1}_{N_i}
			&= \int^1_0 g(r)\big(N_i r-\{N_i r\}+1 - N_i r\big)  \sigma  \,\delta W(r) \nn \\
			&=   \int^1_0  g(r)\big(\frac 12-\{N_i r\} \big) \sigma \, \delta W(r) + \frac 12  \int^1_0 g(r)  \sigma \,  \delta W(r) = Y_{N_i} + \frac 12\int^1_0 g(r)  \sigma\,  \delta W(r).
		\end{align} 
        Substituting \eqref{equ_FubiniMain} into \eqref{equ_exactVn} with $N=N_i$, multiplying both sides by $N_i$, and using
\eqref{equ_Mainpart} yields
\begin{align}\label{equ_decompositionNE}
N_i \widetilde{\mathscr E}_{N_i}(1)&= N_i\uppercase\expandafter{\romannumeral 1}_{N_i}+N_i\uppercase\expandafter{\romannumeral 2}_{N_i}+N_i\uppercase\expandafter{\romannumeral 3}_{N_i}+N_i\uppercase\expandafter{\romannumeral 4}_{N_i} \nn \\
&=  Y_{N_i} +\frac 12\int^1_0 g(r)  \sigma\,  \delta W(r)+N_i\uppercase\expandafter{\romannumeral 2}_{N_i}+N_i\uppercase\expandafter{\romannumeral 3}_{N_i} +N_i\uppercase\expandafter{\romannumeral 4}_{N_i}, 
		\end{align}
   where   we denote the drift term $\uppercase\expandafter{\romannumeral 4}_N:=\int^1_0  g(s) f(X^N(\eta_N(s)))(s - \eta_N(s))  ds$ for notation convenience. 
    The following proposition shows that all terms other than $Y_{N_i}$ on the right-hand side of \eqref{equ_decompositionNE} are negligible under the LIL normalization.
        \begin{prop}\label{prop_Predominantterm}
		 For any
		 strictly increasing	sequence of positive integers $(N_i)_{i\in\N^+}$, we have 
		\begin{align*}
			\lim_{i\to \infty} \frac{ |N_i  \widetilde{\mathscr E}_{N_i} (1)-Y_{N_i} |}{\sqrt{\log\log N_i}} = 0,\quad \text{a.s.}
		\end{align*} 
	\end{prop}
    To prove Proposition~\ref{prop_Predominantterm}, we first
establish the following temporal regularity estimates for $g$
and its Malliavin derivative. The proof is deferred to
Section~\ref{subsec_proofofprop_gDrg}. 
     \begin{prop}\label{prop_gDrg}
            For the function $g$ defined in \eqref{eq:defg} and any $p\geq 2$, there exists a constant $C_p>0$ such that for all $s,t\in[0,1]$,
            \begin{align}
            \Big(  \E \Big[ \big\| g(s)-g(t)\big\|_{\op}^{p} \Big]\Big)^{\frac 1p} &\leq C_p|s-t|^{\frac 12}, \label{equ_Lpofg}\\
                \Big(  \E \Big[ \big(\int^1_0 \big\| \mathcal D_rg(s)-\mathcal D_rg(t)\big\|_{\op}^{2} dr \big)^{\frac p2} \Big]\Big)^{\frac 1p} &\leq C_p|s-t|^{\frac 12} \label{equ_LpofDg}, 
            \end{align}
            where $\mathcal D$ denotes the Malliavin derivative operator. 
        \end{prop}
       
	\begin{proof} 	[Proof of Proposition~\ref{prop_Predominantterm}]
       
	Since $| \int^1_0 g(r) \sigma \, \delta W(r)|<\infty$ almost surely and is independent of $N_i$, we have
		\begin{align}\label{equ_dWO1}
			\lim_{i\to \infty}\frac 12 \frac{|   \int^1_0 g(r)  \sigma  \,\delta W(r)|}{\sqrt{\log \log N_i}} = 0, \quad \text{a.s.}
		\end{align} 
It remains to estimate the other three remainder terms in \eqref{equ_decompositionNE}.

        {\textit{Step 1: Prove $|N\uppercase\expandafter{\romannumeral 2}_N| \to 0$ as $N\to \infty $ almost surely.}} 
        By defining
		$$
		u_N(r)
		:=
		\int_r^{\zeta_N(r)}
		(g(s)-g(r))\sigma ds,
		$$
		we have $
		\uppercase\expandafter{\romannumeral 2}_N=
		\delta(u_N),
		$
		where $\delta$ denotes the divergence operator.
		It follows from the Minkowski integral inequality  and \eqref{equ_Lpofg} that 
		\begin{align*} 
			\big(\E\big[\|u_N(r)\|_{\op}^p\big]\big)^{\frac 1p}
			&\leq \|\sigma\|_{\op }
			\int_r^{\zeta_N(r)}
			\big(\E\big[\|g(s)-g(r)\|_{\op}^p\big]\big)^{\frac 1p} ds\leq
			C_p\int_r^{\zeta_N(r)}
			(s-r)^{1/2} ds
			\leq
			C_pN^{-3/2},
		\end{align*}
		which further gives
		\begin{align}\label{equ_ELpofu}
			\Big(\mathbb E	\Big[	\Big(	\int_0^1\|u_N(r)\|_{\op}^2 dr\Big)^{\frac p2}\Big] \Big)^{1/p}
			\leq
			\Big(
			\int_0^1
			\big(\E \big[ \|u_N(r)\|_{\op}^p \big] \big)^{\frac 2p}\,\mathrm dr
			\Big)^{1/2}\leq
			C_pN^{-3/2}.
		\end{align}
		Using the fact that
		$$
		\mathcal D_\theta u_N(r)
		=
		\int_r^{\zeta_N(r)}
		\bigl(
		\mathcal D_\theta g(s)-\mathcal D_\theta g(r)
		\bigr) \sigma ds,
		$$
		\eqref{equ_LpofDg}, and the Minkowski integral inequality, we derive
		\begin{align*}
			\Big(\mathbb E
			\Big[	\Big(	\int_0^1\|\mathcal D_\theta u_N(r)\|_{\op}^2 d\theta
			\Big)^{\frac p2}	\Big]	\Big)^{\frac 1p}
			&\leq	\|\sigma\| \int_r^{\zeta_N(r)}
			\Big(	\mathbb E
			\Big[	\Big(	\int_0^1	\|\mathcal D_\theta g(s)-\mathcal D_\theta g(r)\|_{\op}^2	d\theta
			\Big)^{\frac p2}	\Big]	\Big)^{\frac 1p}	ds
			\nonumber\\
			&\leq
			C_p
			\int_r^{\zeta_N(r)}
			(s-r)^{\frac 12} ds
			\leq
			C_pN^{-\frac 32},
		\end{align*}
		and then
		\begin{align}	\label{equ_ELpofDu}
			\Big(	\mathbb E
			\Big[	\Big(	\int_0^1\int_0^1
			\| \mathcal D_\theta u_N(r)\|_{\op}^2
			d\theta dr\Big)^{\frac p2}\Big]	\Big)^{\frac 1p}
			&\leq 
			\Big(	\int_0^1
			\Big(\mathbb E
			\Big[	\Big(	\int_0^1 \| \mathcal D_\theta u_N(r)\|_{\op}^2 d\theta
			\Big)^{\frac p2}	\Big]	\Big)^{\frac 2p}
			dr	\Big)^{\frac 12}
			\nonumber\\
			&\leq
			C_pN^{- \frac 32}.
		\end{align}
		Combining \eqref{equ_ELpofu}, \eqref{equ_ELpofDu}, and the Meyer inequality (see e.g. \cite[Theorem 5.3.1]{Nualart2018book})  gives
		\begin{align*}
			\E\big[\big| \uppercase\expandafter{\romannumeral 2}_N \big|^p\big] 
			&\leq
			C_{p}
			\Big(
			\mathbb E	\Big[	\Big(	\int_0^1\|u_N(r)\|_{\op}^2 dr\Big)^{\frac p2}\Big] 
			+
			\mathbb E
			\Big[	\Big(	\int_0^1\int_0^1
			\|\mathcal D_\theta u_N(r)\|_{\op}^2
			d\theta dr\Big)^{\frac p2}\Big]	
			\Big) 		\leq
			C_{p}N^{-\frac {3p}{2}}.
		\end{align*}
		Hence, for any given $\varepsilon>0$,  it follows from the  Chebyshev inequality that
		\begin{align*}
			\sum^\infty_{N=1}	\mathbb P\Big\{ \big|N\uppercase\expandafter{\romannumeral 2}_N \big|>\varepsilon  \Big\}
			&\leq \sum^\infty_{N=1}
			\varepsilon^{-p}
			N^p
			\mathbb E\big[|\uppercase\expandafter{\romannumeral 2}_N |^p\big]
			\leq	C_{p,\varepsilon} \sum^\infty_{N=1} N^{-\frac p2} <\infty,
		\end{align*}
		for $p>2$. Then 
		the Borel--Cantelli lemma implies 
		$$
		\lim_{N\to \infty} | N\uppercase\expandafter{\romannumeral 2}_N|  = 0,
		\quad\text{a.s.}
		$$
		
		{\textit{Step 2: Prove $\lim_{i\to \infty}\frac{| N_i\uppercase\expandafter{\romannumeral 3}_{N_i}| }{\sqrt{\log\log N_i} }= 0$ almost surely.}}
		With the uniform boundedness of $\|\mathcal J(1,z)\|_{\op}, \allowbreak \|Df\|_{\op}, \|D^2f\|_{\op}$ and $\|\sigma\|_{\op}$,  it follows from \eqref{equ_Drg} that
		$
			\| \mathcal D_r g(s) [e_j]\|_{\op} \leq C
		$ with a deterministic constant $C$ independent of $s$, $r$.  Therefore, we have   
		\begin{align*}
			\limsup_{i\to \infty } \frac{| N_i\uppercase\expandafter{\romannumeral 3}_{N_i}| }{\sqrt{\log\log N_i} }  &\leq 	\limsup_{N_i\to \infty } \frac{1}{\sqrt{\log\log N_i}}  \sum^m_{j=1}\| \sigma \|_{\op} N_i  \int^1_0  \int^s_{\eta_{N_i}(s)} \|   \mathcal D_rg(s) [e_j] \|_{\op }   drds \\
			&\leq 	\limsup_{i\to \infty } \frac{C}{\sqrt{\log\log N_i}}  N_i  \int^1_0  (s-{\eta_{N_i}(s)} )ds 
			\leq 	\limsup_{i\to \infty } \frac{C}{\sqrt{\log\log N_i}}  =0, \quad \text{a.s.} 
		\end{align*}

        \textit{Step 3: Prove $\lim_{i \to \infty} \frac{ |N_i \uppercase\expandafter{\romannumeral 4}_{N_i}| }{ \sqrt{\log\log N_i} }=0$ almost surely.}      
		For almost every fixed $\omega$, we first prove that $\sup_{t\in[0,1]} |X^N(t,\omega)|\leq C(\omega)$. 
		It follows from \eqref{equ_EMNon} and the Lipschitz property of $f$ that 
		\begin{align*}
			\sup_{s\in[0,t]} |X^N(s,\omega )| &\leq  |X_0| + C \int^t_0 ( 1+\sup_{r\in[0, s]} |X^N(r,\omega)| ) ds + \|\sigma\|_{\op} \sup_{s\in[0,1]}|W(s,\omega) |\\
			& \leq C(\omega) + C  \int^t_0  \sup_{r\in[0, s]} |X^N(r,\omega)|  ds .
		\end{align*}
		Therefore,  the Gr\"onwall inequality yields $ 	\sup_{s\in[0,1]} |X^N(s,\omega )| \leq C(\omega)$, where the constant $C(\omega)$ is independent of $N$. 

		Recall that  $\|\mathcal J(t,s) \|_{\op} \leq C$ for $0\leq s\leq t\leq 1$. 
		Therefore, 
		\begin{align*}
			 |IV_{N}|&=\Big| 	 \int^1_0g(s) f(X^N(\eta_N(s))) (s-\eta_N(s)) ds \Big| \\
			& \leq C \int^1_0( 1+|X^N(\eta_N(s))|) (s-\eta_N(s))ds\\
			&\leq C(1+\sup_{t\in[0,1]} |X^N(t,\omega)|)\int^1_0(s-\eta_N(s))ds\leq C(\omega)\frac{1}{N},
		\end{align*}
		which leads to
		\begin{align*} 
			\lim_{i \to \infty} \frac{ N_i \big| IV_{N_i}\big|}{ \sqrt{\log\log N_i} }=0, \quad \text{a.s.}
		\end{align*}

		Combining \eqref{equ_dWO1} with {\textit{Step 1--Step 3}}, we have
		\begin{align*}
			\lim_{i \to \infty} \frac{ |N_i \widetilde{\mathscr E}_{N_i}(1) -Y_{N_i}|}{\sqrt{\log\log N_i}} =0,\quad \text{a.s.}, 
		\end{align*}
		which ends the proof.
	\end{proof}
	We conclude this section with the proof of \Cref{prop_remainder}.
	\begin{proof}[Proof of \Cref{prop_remainder}]
			Combining Propositions  \ref{prop_linearization} and \ref{prop_Predominantterm}, we obtain
		\begin{align*}
			&\quad \limsup_{i \to \infty} \frac{ |N_i {\mathcal E}_{N_i}(1) -Y_{N_i}|}{\sqrt{\log\log N_i}} \\
			&\leq 	\limsup_{i \to \infty} \frac{ |N_i ({\mathcal E}_{N_i} (1)-\widetilde{\mathscr E}_{N_i}(1) )|}{\sqrt{\log\log N_i}}  + 	\limsup_{i \to \infty} \frac{ |N_i \widetilde{\mathscr E}_{N_i}(1) -Y_{N_i}|}{\sqrt{\log\log N_i}}  =0, \quad \text{a.s.}, 
		\end{align*} 
		which completes the proof of \Cref{prop_remainder}.
	\end{proof}

	\section{Proof of the almost sure asymptotics of the leading term $Y_{N_i}$}\label{sec_proofofThm}
	In this section, we give the proof of \Cref{prop_LILofY}. 
	For notational simplicity, we denote $G(r):=(U( r))^{-1} Df(X(r))$ in this section.
	Recalling that $\mathcal J(1,r)=U(1)(U(r))^{-1}$, we obtain from the product formula for the divergence operator that 
	\begin{align}\label{equ_exYN}
		Y_N &= \int^1_0U(1) (U( r))^{-1} Df(X(r)) \big(\frac 12 -\{ Nr\}\big) \sigma \,\delta W(r)  \nn\\
		&= U(1)  \int^1_0 G(r) \big(\frac 12 -\{ Nr\}\big) \sigma \,\delta W(r) - \sum^m_{j=1} \int^1_0( \mathcal D_r U(1)[e_j])( G(r)\sigma e_j) \big(\frac 12 -\{ Nr\}\big) dr \nn\\
		&= U(1) 	\mathcal{M}_{N}(1) - \sum^m_{j=1} \int^1_0( \mathcal D_r U(1)[e_j])( G(r)\sigma e_j) \big(\frac 12 -\{ Nr\}\big) dr, 
	\end{align}
	where we define 
	\begin{align*} 
		\mathcal{M}_{N}(t) := \int^t_0  G(r)\big(\frac 12 -\{ Nr\}\big) \sigma dW(r).
	\end{align*}
	For  simplicity, we write $\mathcal{M}_{N}:= 	\mathcal{M}_{N}(1)$. 
	Note that $(\mathcal{M}_{N}(t))_{t\in[0,1]}$ is a  martingale with respect to the filtration ${(\mathcal F_t)}_{t\in[0,1]}$. The decomposition \eqref{equ_exYN} isolates the non-adapted factor $U(1)$ and allows us to analyze the leading term through $\mathcal M_N$.
	We state the almost sure asymptotics for $(\mathcal M_{N_i})_{i\in\N^+}$ in the following proposition. The proofs of the lower and upper bounds for the limit superior are postponed to Sections~\ref{subsec_lowerLIL} and~\ref{subsec_upperLIL}, respectively, while the limit inferior follows by applying the same argument to $-\mathcal M_N$.
		\begin{prop}\label{prop_LILofM}
		Let $v\in\RR^d$ and $(N_i)_{i\in\mathbb N^+}=(\alpha^i)_{i\in\mathbb N^+}$ for a fixed integer \(\alpha\ge2\). Then it holds that
		\begin{align*}
			\limsup_{i\to \infty} \frac{ v^\top \mathcal M_{N_i}}{\sqrt{\log\log N_i} }=\sqrt{2v^\top \mathcal{V}_{M,\sigma} v},\quad \text{a.s.}\quad \text{and} \quad \liminf_{i\to \infty} \frac{ v^\top \mathcal M_{N_i}}{\sqrt{\log\log N_i} }=-\sqrt{2v^\top \mathcal{V}_{M,\sigma} v},\quad \text{a.s.},
		\end{align*}
		where the asymptotic covariance matrix is defined by 
		\begin{align}\label{equ_DefVM}
		\mathcal V_{M,\sigma}:= \frac{1}{12}\int^1_0 G(r) \sigma \sigma^\top (G(r))^\top   dr.
		\end{align}
	\end{prop}
	
	With \Cref{prop_LILofM}, we now present the proof of \Cref{prop_LILofY}.
	\begin{proof}[Proof of \Cref{prop_LILofY}]
		We first show that the second term on the right-hand side of \eqref{equ_exYN} is negligible under the LIL normalization. 
For any  fixed $t\in[0,1]$, it follows from \eqref{equ_UEqu} that
		\begin{align*}
			 U (t) = I_d + \int^t_0 Df(X(z))U(z) dz.
		\end{align*}
Differentiating this identity in the Malliavin sense yields, for any $u\in\RR^m$ and  $v\in\RR^d$, 
		\begin{align}\label{equ_DrU}
			(\mathcal D_r U(t) [u])v &=\int^t_0   \Big[\big(\mathcal D_r \big(Df(X(z))\big) [u] \big)(U(z)v) + Df(X(z)) (\mathcal D_rU(z)  [u])v\Big] dz.
		\end{align}
        Thus, $ (\mathcal D_r U(t) [u])v$ satisfies
		\begin{align*}
			\frac{d}{dt} (\mathcal D_r U(t) [u])v&= (\mathcal D_r(Df(X(t))) [u]) (U(t)v)+ Df(X(t)) (\mathcal D_r U(t)[u])v, \\
			(\mathcal D_r U(r) [u])v &=\mathbf{0}_d.
		\end{align*}
        Here, $u$  specifies the noise direction, and \(v\) is the vector on which the resulting matrix acts.
		Combining this with the expression for $\mathcal D_r X$ in \eqref{equ_DrX},  for $0\leq r\leq t$, $u\in\RR^m, v\in\RR^d$, we obtain
		\begin{align*}
			(\mathcal D_r U(t) [u])v &= \int^t_0 \mathcal J(t,z)  (\mathcal D_r(Df(X(z))) [u]) (U(z)v) dz \\
			&=\int^t_r  \mathcal J(t,z) D^2f(X(z))(\mathcal J(z,r)\sigma u, U(z)v) dz.
		\end{align*}
		Therefore, the uniform boundedness of the operators leads to 
		\begin{align}\label{equ_remainerY}
			\lim_{i\to\infty} \frac{1}{\sqrt{\log\log N_i}}\sum^m_{j=1} \Big|   \int^1_0( \mathcal D_r U(1)[e_j])( G(r) \sigma e_j) \big(\frac 12 -\{ N_ir\}\big) dr\Big| =0 ,\quad \text{a.s.}
		\end{align}
	Combining \eqref{equ_exYN} and \eqref{equ_remainerY}, we arrive at
		\begin{align}\label{equ_LILYtoM}
		\limsup_{i\to\infty} \frac{   v^\top Y_{N_i} }{\sqrt{\log\log N_i} }=\limsup_{i\to\infty} \frac{   ((U(1))^\top v)^\top \mathcal{M}_{N_i} }{\sqrt{\log\log N_i} }, \quad \text{a.s.}
	\end{align}

Since $\mathbb{Q}^d$ is a countable dense subset of $\RR^d$, it follows from \Cref{prop_LILofM} that there exists $\Omega_0 \subset \Omega$ with $\PP(\Omega_0)=1$ such that for any $\omega\in\Omega_0$, 
\begin{align}\label{equ_LILQd}
	\limsup_{i\to\infty} \frac{  \xi^\top \mathcal{M}_{N_i}(\omega) }{\sqrt{\log\log N_i} }  = \sqrt{2\xi^\top \mathcal{V}_{M,\sigma}(\omega) \xi}, \quad 	\liminf_{i\to\infty} \frac{  \xi^\top \mathcal{M}_{N_i}(\omega) }{\sqrt{\log\log N_i} }  = -\sqrt{2\xi^\top \mathcal{V}_{M,\sigma}(\omega) \xi}, \qquad \forall \xi\in \mathbb{Q}^d.
\end{align}
	Fix $\omega\in\Omega_0$.
Applying \eqref{equ_LILQd} to the coordinate vectors yields
 $ 	\limsup_{i\to\infty} \frac{  | \mathcal{M}_{N_i}(\omega) |}{\sqrt{\log\log N_i} } \leq C(\omega)$.
Taking $\xi_1,\xi_2\in\RR^d$, it holds that
\begin{align*}
	\Big|	\limsup_{i\to\infty} \frac{  \xi_1^\top \mathcal{M}_{N_i}(\omega) }{\sqrt{\log\log N_i} }  - \limsup_{i\to\infty} \frac{  \xi_2^\top \mathcal{M}_{N_i}(\omega) }{\sqrt{\log\log N_i} }\Big|  \leq \limsup_{i\to\infty} \Big| \frac{  (\xi_1-\xi_2)^\top \mathcal{M}_{N_i}(\omega) }{\sqrt{\log\log N_i} }  \Big| \leq C(\omega) | \xi_1-\xi_2| ,
\end{align*}
	which shows the continuity of the mapping $\xi\mapsto \limsup_{i\to\infty} \frac{   \xi^\top \mathcal{M}_{N_i}(\omega) }{\sqrt{\log\log N_i} } $.   Combining this with the continuity of the mapping $\xi\mapsto  \sqrt{2\xi^\top \mathcal{V}_{M,\sigma}(\omega) \xi}$, we extend \eqref{equ_LILQd} to all $\xi\in\RR^d$. Specifically, for any $\omega\in\Omega_0$, we have
	\begin{align*}
		\limsup_{i\to\infty} \frac{  \xi^\top \mathcal{M}_{N_i}(\omega) }{\sqrt{\log\log N_i} }  = \sqrt{2\xi^\top \mathcal{V}_{M,\sigma}(\omega) \xi}, \quad 	\liminf_{i\to\infty} \frac{  \xi^\top \mathcal{M}_{N_i}(\omega) }{\sqrt{\log\log N_i} }  = -\sqrt{2\xi^\top \mathcal{V}_{M,\sigma}(\omega) \xi}, \qquad \forall \xi\in \mathbb{R}^d.
	\end{align*}
	For every $\omega\in\Omega_0$ and $v\in\mathbb R^d$,
	we take $\xi=(U(1,\omega))^\top v$ and obtain
	\begin{align}\label{equ_omegaLIL}
		\limsup_{i\to\infty}
		\frac{ ((U(1,\omega))^\top v)^\top\mathcal M_{N_i}(\omega)}
		{\sqrt{\log\log N_i}}
		&=
		\sqrt{
			2v^\top U(1,\omega)
			\mathcal V_{M,\sigma}(\omega)
			U(1,\omega)^\top v
		}=
		\sqrt{2v^\top\mathcal V_\sigma(\omega)v},
	\end{align}
	where we use the identity
	$
	U(1)\mathcal V_{M,\sigma}U(1)^\top
	=\mathcal V_\sigma.
	$
	By plugging \eqref{equ_omegaLIL} into \eqref{equ_LILYtoM}, we obtain the limit superior identity in \Cref{prop_LILofY}. The proof is complete, since 
	the corresponding limit inferior identity follows similarly. 
	\end{proof}
	\subsection{Lower bound estimate in \Cref{prop_LILofM}}\label{subsec_lowerLIL}
	We now prove the lower bound for the limit superior
	in \Cref{prop_LILofM}.
	We divide the proof into three steps. First, we construct a conditionally Gaussian sequence ${(\mathcal M_i^{(l)})}_{i\geq l}$  and establish exponential decay of its conditional covariances. Second, we apply a conditional version of the Kochen--Stone lemma to derive the LIL lower bound for this sequence. Finally, 
    we estimate the difference between \(\mathcal M_i^{(l)}\) and \(\mathcal M_{N_i}\) to obtain the desired lower bound for \((\mathcal M_{N_i})_{i\in\mathbb N^+}\).
	\subsubsection{Construction of the conditionally Gaussian sequence $(\mathcal M_i^{(l)})_{i\geq l}$}
	We first present the construction and give some basic estimates.
	We denote by $X^{(l)}_k$ the numerical solution defined by \eqref{equ_EMNon} at $t^{(l)}_k=kh_l$ with step size $h_l:=N_l^{-1}=\alpha^{-l}$, where  $l\in\N^+$. We write $\Delta W^{(l)}_k:=W(t^{(l)}_{k+1})-W(t^{(l)}_k)$ and define
	$
	\mathcal G^{(l)}:=\sigma\{\Delta W^{(l)}_k:0\leq k<N_l\}.
	$
	For simplicity, we write $t_k=t^{(l)}_k$ when no confusion arises.  
	Define $$\mathcal M^{(l)}_i(t) := \int^t_0 G^{(l)}(r) \big( \frac 12 - \{ N_i r\}\big)\sigma d W(r),\quad t\in [0,1]$$
	and write $\mathcal M^{(l)}_i:=\mathcal M^{(l)}_i(1)$ for simplicity, where
	 $G^{(l)}(r)=G^{(l)}(t_k)$ for $r\in[t_k, t_{k+1})$ with 
	\begin{align*}
		G^{(l)}(t_k):=(U^{(l)}(t_k))^{-1}Df(X^{(l)}_k).
	\end{align*}
	Here, $U^{(l)}(\cdot)$ is the solution of
	\begin{align}\label{equ_discreteU}
		dU^{(l)}(r) = Df(X^{(l)}(\eta_{N_l}(r))) U^{(l)}(r) dr,\quad  U^{(l)}  (0)=I_d.
	\end{align} 
	In the following, we show the uniform boundedness and temporal continuity of $U^{(l)}$.
	By solving \eqref{equ_discreteU} explicitly, we obtain
	\begin{align*}
		U^{(l)} (t_{k+1}) = e^{h_l Df(X^{(l)}_k) } U^{(l)} (t_{k}) = \cdots =e^{h_l Df(X^{(l)}_k) }e^{h_l Df(X^{(l)}_{k-1}) }\cdots e^{h_l  Df(X^{(l)}_0)}.
	\end{align*}
	Then
	the uniform boundedness of $U^{(l)}$ follows from the uniform  boundedness of $\|Df\|_{\op}$ as follows 
	\begin{align}\label{equ_boundedUl}
		\sup_{l\in\N^+}\sup_{t\in[0,1]}\| U^{(l)}(t)\|_{\op}  
        \leq \sup_{l\in \N^+} \prod_{j=0}^{N_l-1} e^{h_l  \| Df(X^{(l)}_j)\|_{\op} } 
        \leq C. 
	\end{align}
	The temporal continuity follows from the fact that 
	$  \sup_{t\in[0,1]}  \| U^{(l)}(t)-U^{(l)}(\eta_{N_l} (t)) \|_{\op} \leq C\frac{1}{N_l}.$
	Applying the same arguments as in the estimate of $U^{-1}$ yields the uniform boundedness of $(U^{(l)})^{-1}$ and the continuity estimate
	\begin{align}\label{equ_continuityUlInverse}
		\sup_{t\in[0,1]}  \| (U^{(l)}(t))^{-1}-(U^{(l)}(\eta_{N_l} (t)) )^{-1}\|_{\op} \leq C\frac{1}{N_l}.
	\end{align}
	
	Next we show that $(\mathcal M^{(l)}_i)_{i\geq l}$ is a conditionally Gaussian sequence and give its conditional covariance. 
   Any finite collection of the stochastic integrals
$\int_{t_k}^{t_{k+1}}(\frac12-\{N_i r\})\sigma dW(r),~i\geq l,\, 0\leq k\leq N_l-1$
and the increments generating $\mathcal G^{(l)}$ are jointly
Gaussian. Moreover, for $k\in\{0,\ldots,N_l-1\}$ and $i\geq l$, 
we have
	\begin{align*}
		\text{Cov}\Big( \int^{t_{k+1}}_{t_k} \big(\frac 12-\{ N_i r\} \big) \sigma dW(r), \Delta W^{(l)}_k\Big) &=	\E \Big[ \big( \int^{t_{k+1}}_{t_k} \big(\frac 12-\{ N_i r\} \big) \sigma dW(r) \big)  \big(\int^{t_{k+1}}_{t_k}  dW(r) \big)^\top\Big] \\
		&= \sigma  \int^{t_{k+1}}_{t_k} \big(\frac 12-\{ N_i r\} \big) dr =\mathbf 0_{d\times m}.
	\end{align*}
The covariances with the remaining coarse-grid increments also vanish, by the independence of Wiener increments over disjoint intervals. Therefore, the family of the stochastic integrals is independent of $\mathcal G^{(l)}$.
	Then, for $i,j\geq l$, we have 
	\begin{align*}
	&\quad	\E\big[ \mathcal M^{(l)}_i (\mathcal M^{(l)}_j)^\top \, \big| \, \mathcal G^{(l)} \big]\\
		&=	\E\Big[ \Big( \sum^{N_l-1}_{k_1=0}G^{(l)}(t_{k_1}) \int^{t_{{k_1}+1}}_{t_{k_1}}   \big( \frac 12 - \{ N_i r\}\big)\sigma dW(r) \Big)\Big( \sum^{N_l-1}_{k_2=0}G^{(l)}(t_{k_2}) \int^{t_{{k_2}+1}}_{t_{k_2}}   \big( \frac 12 - \{ N_j r\}\big)\sigma d W(r)\Big)^\top \, \Big| \, \mathcal G^{(l)} \Big] 
		\\
		&=	\sum^{N_l-1}_{k_1,k_2=0} G^{(l)}(t_{k_1})  \E\Big[ \Big(  \int^{t_{{k_1}+1}}_{t_{k_1}}  \big( \frac 12 - \{ N_i r\}\big)\sigma dW(r) \Big)\Big( \int^{t_{{k_2}+1}}_{t_{k_2}}   \big( \frac 12 - \{ N_j r\}\big)\sigma dW(r)\Big)^\top  \Big] (G^{(l)}(t_{k_2}))^\top 
		\\
		&=\sum^{N_l-1}_{k=0} (G^{(l)}(t_k))\sigma \sigma^\top (G^{(l)}(t_k))^\top \int^{t_{k+1}}_{t_k}   \big(\frac 12-\{ N_i r\} \big)\big(\frac 12-\{ N_j r\} \big) dr.
	\end{align*}
	Since $ \frac 12 - \{x\}$ is a periodic function with $T=1$,  the Fourier expansion gives
	\begin{align*} 
		\frac 12 - \{ x \}= \sum_{m=1}^\infty \frac{\sin(2\pi m x)}{\pi m}
	\end{align*}
 for all $x\notin \mathbb Z$, 
	which, combined with the identity  $2\sin x\sin y=\cos(x-y)-\cos(x+y)$, leads to
	\begin{align*}
		\big(\frac 12-\{ N_i r\} \big)\big(\frac 12-\{ N_j r\} \big)&= \sum_{m=1}^\infty\sum_{n=1}^\infty \frac{\sin(2\pi m N_i r)}{\pi m}  \frac{\sin(2\pi n N_jr )}{\pi n} \nn \\
		&= \frac{1}{2\pi^2} \sum_{m=1}^\infty\sum_{n=1}^\infty  \frac{1}{mn}\big( \cos (2\pi(mN_i-nN_j)r) -  \cos (2\pi(mN_i+nN_j)r) \big).
	\end{align*}
	Therefore, we have
	\begin{align*}
		&\quad   \int^{t_{k+1}}_{t_k}   \big(\frac 12-\{ N_i r\} \big)\big(\frac 12-\{ N_j r\} \big) dr\\
		&=\frac{1}{2\pi^2} \sum_{m=1}^\infty\sum_{n=1}^\infty  \frac{1}{mn}\int^{t_{k+1}}_{t_k} \big( \cos (2\pi(mN_i-nN_j)r) -  \cos (2\pi(mN_i+nN_j)r) \big) dr. 
	\end{align*}
	For the second part, the condition $i,j\geq l$ gives
	\begin{align} \label{equ_int2ndterm_Ik}
		&\quad \int^{t_{k+1}}_{t_k}  \cos (2\pi(mN_i+nN_j)r)   dr \nn \\
		&= \frac{1}{2\pi (mN_i+nN_j)}\big( \sin(2\pi(mN_i+nN_j) t_{k+1}) - \sin(2\pi(mN_i+nN_j) t_{k}) \big) \nn \\
		&= \frac{1}{2\pi (mN_i+nN_j)}\big( \sin(2\pi (k+1)(m\alpha^{i-l}+n\alpha^{j-l})) - \sin(2\pi k(m\alpha^{i-l}+n\alpha^{j-l}))=0.
	\end{align}
	When $mN_i\neq nN_j$,  applying an  argument similar to that used in \eqref{equ_int2ndterm_Ik}, we derive $$\int^{t_{k+1}}_{t_k}  \cos (2\pi(mN_i-nN_j)r)   dr=0.$$ When $mN_i =  nN_j$, it holds that
	\begin{align*}
		\int^{t_{k+1}}_{t_k}  \cos (2\pi(mN_i-nN_j)r)   dr = h_l.
	\end{align*}
	
	Combining results above gives
	\begin{align*}
		\int^{t_{k+1}}_{t_k}   \big(\frac 12-\{ N_i r\} \big)\big(\frac 12-\{ N_j r\} \big) dr = \frac{h_l}{2\pi^2} \sum_{\substack{m,n=1  \\ mN_i =  nN_j  }}^\infty  \frac{1}{mn} = \frac{h_l}{2\pi^2}\sum_{n=1}^\infty \frac{1}{n^2 \alpha^{|j-i| }} = \frac{h_l}{12}\alpha^{-|j-i|}.
	\end{align*}
	This therefore yields
	\begin{align*}
		\E\big[ \mathcal M^{(l)}_i (\mathcal M^{(l)}_j)^\top \, \big| \, \mathcal G^{(l)} \big] &= \alpha^{-|j-i|}  \frac{h_l}{12} \sum^{N_l-1}_{k=0} (G^{(l)}(t_k))\sigma \sigma^\top (G^{(l)}(t_k))^\top \\
		&=  \alpha^{-|j-i|} \frac{1}{12} \int^1_0  (G^{(l)}(r))\sigma \sigma^\top(G^{(l)}(r))^\top  dr 
		=  \alpha^{-|j-i|}  \mathcal V^{(l)}_{M,\sigma},
	\end{align*}
	where we define 
	\begin{align}\label{equ_DefVMl}
		\mathcal V^{(l)}_{M,\sigma}:= \frac{1}{12} \int^1_0  (G^{(l)}(r))\sigma \sigma^\top(G^{(l)}(r))^\top  dr.
	\end{align}
	Note that the factor $\alpha^{-|j-i|}$ here yields exponential decay of the	covariance as the index gap $|j-i|$ increases, which is key to
	verifying the weak dependence condition in the following proof.
	\subsubsection{Lower bound  for the LIL of $(\mathcal M^{(l)}_i)_{i\geq l}$}
On the event \((v^\top\mathcal V^{(l)}_{M,\sigma}v=0)\), we have \(v^\top\mathcal M_i^{(l)}=0\) almost surely for all \(i\ge l\), so the desired lower bound holds trivially. We therefore work on the event \((v^\top\mathcal V^{(l)}_{M,\sigma}v>0)\) throughout the following argument.
	We first state a conditional version of the Kochen--Stone lemma \cite{KochenStone1964};
	see, e.g., \cite[Lemma 5.3]{ConditionalKSthm2022}.
	\begin{lemma}[Conditional version of Kochen--Stone lemma] \label{lemma_Kochen--Stone}
		Let  ${(A_n)}$ be a sequence of events on a probability space $(\Omega,\mathcal F,\PP)$ and  $\mathcal G\subseteq \mathcal F$ a $\sigma$-algebra.
		Then, 	if
		$\sum_{n=1}^\infty \PP\{A_n\,  |\,  \mathcal G\} = \infty$ almost surely, we have
		\begin{align} \label{equ_Kochen--Stone}
		\PP\Big\{A_n \;\; \text{i.o.}\, \Big| \,\mathcal G\Big\} \geq	\limsup_{k\to\infty} \frac{ \big(\sum_{n=1}^k \PP\{A_n\, | \, \mathcal G\}\big)^2 }{\sum_{1\le m,n\le k} \PP\{A_m\cap A_n \,|\, \mathcal G\} },\qquad\text{a.s.}
		\end{align}
	\end{lemma}
    The Kochen--Stone lemma extends the second Borel--Cantelli lemma
by providing a lower bound for the probability that infinitely many events occur without requiring independence. This allows us to establish the almost sure result under a weak correlation condition.
For any	 fixed $v\in\RR^d$ and  $\varepsilon\in(0,1)$, we define  $m^{(l,v)}:= (1-\varepsilon) \sqrt{2 v^\top \mathcal V^{(l)}_{M,\sigma} v} $ and 
	$	c^{(l,v)}_{N_i}:= m^{(l,v)}\sqrt{\log\log N_i}$,
	which are $\mathcal G^{(l)}$-measurable random variables. 
	With $(v^\top \mathcal M^{(l)}_{i}> c^{(l,v)}_{N_i})$ and $\mathcal G^{(l)}$ serving as $(A_n)$ and $\mathcal G$ in \Cref{lemma_Kochen--Stone}, in the following we establish the divergence of the series $\sum_{n=1}^\infty \mathbb{P}\{A_n \mid \mathcal G\} = \infty$ in (i) and prove that $\PP\big\{A_n \;\; \text{i.o.}\, \big| \,\mathcal G\big\} =1$ via \eqref{equ_Kochen--Stone} in (ii), respectively. 
	
	(\romannumeral 1) 
	We first estimate  $\sum^\infty_{i=i_{0}}\PP\{ v^\top \mathcal M^{(l)}_{i}> c^{(l,v)}_{N_i} \, | \, \mathcal G^{(l)}\}$ for some $i_{0}\geq l$ large enough such that $\log\log N_i\geq 1$ for any $i\geq i_{0}$.  Note that $v^\top \mathcal M^{(l)}_{i}$ is a central Gaussian random variable conditionally on $\mathcal G^{(l)}$ with the conditional variance 
	\begin{align*}
		\E[ (v^\top \mathcal M^{(l)}_{i} )^2 | \mathcal G^{(l)}]  =v^\top   \E[  \mathcal M^{(l)}_{i} (\mathcal M^{(l)}_{i} )^\top| \mathcal G^{(l)}]  v = v^\top \mathcal V^{(l)}_{M,\sigma} v. 
	\end{align*}
	An elementary bound for the Gaussian tail probability gives
	$1-\Phi(x) \geq \frac{x}{1+x^2} \frac{e^{-x^2/2}}{\sqrt{2\pi}}\geq  \frac{e^{-x^2/2}}{2\sqrt{2\pi}x}$ for $x\geq 1$, where $\Phi(x)$ denotes the distribution function of a standard Gaussian random variable. For any $\varepsilon\in(0,1-1/\sqrt{2})$, this yields
    \begin{align}\label{equ_seriesconditionalP}
		\sum^\infty_{i=i_0}\PP\Big\{ v^\top \mathcal M^{(l)}_{i}> c^{(l,v)}_{N_i} \, \Big| \, \mathcal G^{(l)} \Big\} 
		&= \sum^\infty_{i=i_0}\PP\bigg\{ \frac{v^\top \mathcal M^{(l)}_{i}}{\sqrt{v^\top\mathcal V^{(l)}_{M,\sigma} v}}> \frac{c^{(l,v)}_{N_i}}{\sqrt{v^\top\mathcal V^{(l)}_{M,\sigma} v}} \, \bigg| \, \mathcal G^{(l)}\bigg\} \nn\\
        &\geq  \sum^\infty_{i=i_0}  \frac{1 }{ 2\sqrt{2\pi}(1-\varepsilon) \sqrt{2\log\log {N_i}} }\exp\Big\{- \frac{(1-\varepsilon)^22\log\log {N_i}}{2} \Big\} \nn \\ 
        &=\sum^\infty_{i=i_0}  \frac{1 }{ 4\sqrt{\pi}(1-\varepsilon) \sqrt{\log i+\log\log \alpha }} (\log \alpha)^{-(1-\varepsilon)^2}  i^{-(1-\varepsilon)^2}  \nn \\ 
		& \geq C({\alpha} )\sum^\infty_{i=i_0} \frac{1}{\sqrt{\log i}} i^{-(1-\varepsilon)^2}=\infty,
	\end{align}
	which leads to the divergence of the series. 

	(\romannumeral 2) 
	Denote $E^{(l,v)}_i:= \big(  v^\top \mathcal M^{(l)}_{i}>c^{(l,v)}_{N_i} \big)$. 
	Then we show that the right-hand side of \eqref{equ_Kochen--Stone} equals one. 
	Following from the conditional Gaussian property of $\mathcal M^{(l)}_i$,  we have the decomposition for $i\neq j$,
	\begin{align*}
		\PP\big\{ E^{(l,v)}_i \cap E^{(l,v)}_j\, \big| \, \mathcal G^{(l)}\big\} = \PP\big\{ E^{(l,v)}_i\, \big| \, \mathcal G^{(l)}\big\} \PP\big\{E^{(l,v)}_j\, \big| \, \mathcal G^{(l)}\big\}  +  \mathcal 	R^{(l,v)}(i,j),
	\end{align*}
	where the correlation coefficient 
	$\rho_{i,j} := \E \big[ \frac{ (v^\top \mathcal M^{(l)}_{i})( v^\top \mathcal M^{(l)}_{j})}{v^\top \mathcal V^{(l)}_{M,\sigma} v }  \big| \mathcal G^{(l)} \big]=\alpha^{-|j-i|} $
	and  we denote 
	$
	\mathcal 	R^{(l,v)}(i,j):=	\int^{\rho_{i,j}}_0 \phi^{(l,v)}(c^{(l,v)}_{N_i},c^{(l,v)}_{N_j},\lambda) d\lambda
	$
	with
	$$\phi^{(l,v)}(c^{(l,v)}_{N_i}, c^{(l,v)}_{N_j}; \lambda) := \frac{1}{2\pi \sqrt{1-\lambda^2}} \exp\bigg\{ -\frac{(c^{(l,v)}_{N_i})^2 + (c^{(l,v)}_{N_j})^2 - 2\lambda c^{(l,v)}_{N_i} c^{(l,v)}_{N_j}}{2(1-\lambda^2)v^\top \mathcal V^{(l)}_{M,\sigma} v } \bigg\}.$$
	Therefore, due to the equivalence expression
	\begin{align*}
		&\quad \liminf_{k\to\infty} \frac{\sum_{1\le i,j\le k} 	\PP\big\{ E^{(l,v)}_i \cap E^{(l,v)}_j\, \big| \, \mathcal G^{(l)}\big\} }{\big(\sum_{i=1}^k \PP\{E^{(l,v)}_i \, | \, \mathcal G^{(l)}\}\big)^2}\\ &= 	\liminf_{k\to\infty} \frac{ 	\big(\sum_{i=1}^k \PP\{E^{(l,v)}_i \, | \, \mathcal G^{(l)}\}\big)^2 +\sum^k_{i=1} \PP\{E^{(l,v)}_i \, | \, \mathcal G^{(l)}\}\big(1-\PP\{E^{(l,v)}_i \, | \, \mathcal G^{(l)}\}\big)+ 2\sum_{1\leq i<j\leq k} \mathcal R^{(l,v)}(i,j) }{\big(\sum_{i=1}^k \PP\{E^{(l,v)}_i \, | \, \mathcal G^{(l)}\}\big)^2 }
	\end{align*}
	and the convergence result for the second term
	\begin{align*}
		\lim_{k\to\infty}\frac{	\sum^k_{i=1} \PP\{E^{(l,v)}_i \, | \, \mathcal G^{(l)}\}\big(1-\PP\{E^{(l,v)}_i \, | \, \mathcal G^{(l)}\}\big) }{ \big(\sum_{i=1}^k \PP\{E^{(l,v)}_i \, | \, \mathcal G^{(l)}\}\big)^2}\leq\lim_{k\to\infty} \frac{1}{\sum_{i=1}^k \PP\{E^{(l,v)}_i \, | \, \mathcal G^{(l)}\}}  =0,
	\end{align*}
	once we prove that
	\begin{align}\label{equ_Rlim}
		\lim_{k\to\infty} \frac{2\sum_{1\leq i<j\leq k} \mathcal R^{(l,v)}(i,j) }{\big(\sum_{i=1}^k \PP\{E^{(l,v)}_i \, | \, \mathcal G^{(l)}\}\big)^2}=0,
	\end{align}
	we arrive at
	\begin{align*}
		\limsup_{k\to\infty} \frac{	\big(\sum_{i=1}^k \PP\{E^{(l,v)}_i \, | \, \mathcal G^{(l)}\}\big)^2 }{\sum_{1\le i,j\le k} 	\PP\big\{ E^{(l,v)}_i \cap E^{(l,v)}_j\, \big| \, \mathcal G^{(l)}\big\}
		}=	\bigg(\liminf_{k\to\infty} \frac{ \sum_{1\le i,j\le k} 	\PP\big\{ E^{(l,v)}_i \cap E^{(l,v)}_j\, \big| \, \mathcal G^{(l)}\big\}}{	\big(\sum_{i=1}^k \PP\{E^{(l,v)}_i \, | \, \mathcal G^{(l)}\}\big)^2
		}\bigg)^{-1}=1.
	\end{align*}
	
	Next we give the proof of \eqref{equ_Rlim}.
	Using
		\begin{align*}
	&\quad \frac{(c^{(l,v)}_{N_i})^2 + (c^{(l,v)}_{N_j})^2 - 2\lambda c^{(l,v)}_{N_i} c^{(l,v)}_{N_j}}{2(1-\lambda^2)}\\
	& \ge \frac{((c^{(l,v)}_{N_i})^2 + (c^{(l,v)}_{N_j})^2)(1-\lambda)}{2(1-\lambda)(1+\lambda)} = \frac{(c^{(l,v)}_{N_i})^2 +( c^{(l,v)}_{N_j})^2}{2(1+\lambda)}\geq \frac{(c^{(l,v)}_{N_i})^2 + (c^{(l,v)}_{N_j})^2}{2(1+\rho_{i,j})},
\end{align*}
	and
	\begin{align*}
		\int^{\rho_{i,j}}_0 \frac{1}{\sqrt{1-\lambda^2 }} d\lambda =\arcsin(\rho_{i,j}) \leq \frac{\pi}{2} \rho_{i,j}
	\end{align*}
	for $\rho_{i,j}\in(0,1)$, 
	we have
	\begin{align*}
		\mathcal R^{(l,v)}(i,j)\leq \frac 14 \rho_{i,j}\exp\Big\{ - \frac{(c^{(l,v)}_{N_i})^2 + (c^{(l,v)}_{N_j})^2}{2(1+\rho_{i,j}) v^\top \mathcal V^{(l)}_{M,\sigma} v} \Big\}.
	\end{align*}
	Since $0\leq \rho_{i,j}\leq \alpha^{-1}<1$, it holds that
	\begin{align*}
		\mathcal 	R^{(l,v)}(i,j) &\leq C \alpha^{-|j-i|} 	\exp\Big\{ -  \frac{(m^{(l,v)})^2(\log i+ \log j+2\log\log \alpha)}{2(1+\alpha^{-1}) v^\top \mathcal V^{(l)}_{M,\sigma} v} \Big\}= C(\alpha) \alpha^{-|j-i|} (ij)^{-\frac{\alpha(1-\varepsilon)^2 }{\alpha+1}  }.
	\end{align*}
	Denote  $\gamma:=
     \frac{\alpha(1-\varepsilon)^2 }{\alpha+1} $.  Since
	\begin{align*}
		\sum_{i_0\leq i<j\leq k} 	\mathcal R^{(l,v)}(i,j) \leq  C\sum^{k}_{i=i_0} \sum^{k}_{j=i+1} \alpha^{-(j-i)} (ij)^{-\gamma }\leq C \sum^{k}_{i=i_0} i^{-2\gamma} \big(\sum^\infty_{j=i+1} \alpha^{-(j-i)}  \big) \leq C \sum^{k}_{i=i_0} i^{-2\gamma},
	\end{align*}
	the partial sums are \(O(1)\) when \(2\gamma>1\), \(O(\log k)\) when \(2\gamma=1\), and \(O(k^{1-2\gamma})\) when \(2\gamma<1\).
As a byproduct of \eqref{equ_seriesconditionalP}, for any integer $k\geq 2i_0$, we have
	\begin{align} \label{equ_seriestoK}
		&\quad \sum^k_{i=i_0}\PP\{ v^\top \mathcal M^{(l)}_{i}> c^{(l,v)}_{N_i} \, | \, \mathcal G^{(l)}\} \nn \\&\geq	C\sum^k_{i=i_0} \frac{1}{\sqrt{\log i}} i^{-(1-\varepsilon)^2}\geq C	\sum^k_{i=[k/2]} \frac{1}{\sqrt{\log i}} i^{-(1-\varepsilon)^2} \geq C \frac{1}{\sqrt{\log k}} k^{1-(1-\varepsilon)^2}.
	\end{align}
	By \eqref{equ_seriestoK}, \eqref{equ_Rlim} follows immediately
when $2\gamma\geq 1$.
	For the case $2\gamma<1$, it follows from \eqref{equ_seriestoK} and the estimate 
	\begin{align*}
		-1+ 2(1-\varepsilon)^2 - 2\gamma = -1 + \frac{2(1-\varepsilon)^2 }{\alpha+1} < 0 
	\end{align*}
	that
	\begin{align*}
		\lim_{k\to\infty} \frac{2\sum_{1\leq i<j\leq k} \mathcal R^{(l,v)}(i,j) }{\big(\sum_{i=1}^k \PP\{E^{(l,v)}_i \, | \, \mathcal G^{(l)}\}\big)^2} \leq 	C \lim_{k\to\infty} (\log k)\, k^{-1+ 2(1-\varepsilon)^2-2\gamma} =0.
	\end{align*} 
	
	Combining (\romannumeral 1)--(\romannumeral 2), applying \Cref{lemma_Kochen--Stone} gives
	$\PP\{v^\top \mathcal M^{(l)}_{i}> c^{(l,v)}_{N_i} \  \text{i.o.} \, | \, \mathcal G^{(l)}  \}=1$. Let $\varepsilon\to0$, we arrive at
	\begin{align}\label{equ_limsupMl}
		\limsup_{i\to \infty} \frac{v^\top \mathcal M^{(l)}_{i}}{\sqrt{\log\log N_i} } \geq \sqrt{2 v^\top \mathcal V^{(l)}_{M,\sigma} v}, \quad \text{a.s.}
	\end{align}
	\subsubsection{Lower bound for the LIL of the original sequence $(\mathcal M_{N_i})_{i\in\N^+}$}
	We first give some useful estimates for $U$, $U^{(l)}$, $G$, and $G^{(l)}$. 
	Combining the uniform mean-square
	convergence estimate \eqref{equ_UniformLpE}
	 with $f\in\mathcal C^\infty_b$ and the H\"older continuity of $X$ in  \eqref{equ_XMomentcontinuity}, we obtain
	\begin{align}\label{equ_Dfappro}
		&\quad \sup_{0\leq r\leq 1}	\E\big[ \|Df(X(r)) -Df(X^{(l)}(\eta_{N_l}(r))) \|_{\op}^2 \big]  \nn \\
		&\leq 2\sup_{0\leq r\leq 1}\Big(	\E\big[ \|Df(X(r)) -Df(X(\eta_{N_l}(r))) \|_{\op}^2 \big] + \E\big[ \|Df(X(\eta_{N_l}(r))) - Df(X^{(l)}(\eta_{N_l}(r)) )\|_{\op}^2 \big] \Big)\nn \\
		&\leq C \sup_{0\leq r\leq 1} (r-\eta_{N_l}(r))+ C \frac{1}{(N_l)^2}\leq C  \frac{1}{N_l}.
	\end{align}
	Based on the equations for $ U$ and $U^{(l)}$ (\eqref{equ_UEqu} and \eqref{equ_discreteU}), we next give the corresponding approximation estimate. It follows from the expression 
	\begin{align*}
		U(t) -U^{(l)}(t) = \int^t_0 Df(X(r))\big( U(r) -U^{(l)}(r)\big) dr +  \int^t_0 \big(Df(X(r))-Df(X^{(l)}(\eta_{N_l}(r)))\big) U^{(l)}(r) dr, 
	\end{align*}
	the uniform boundedness of $Df$ and $U^{(l)}$ \eqref{equ_boundedUl}, and the approximation property \eqref{equ_Dfappro} that for any $t\in[0,1]$,
	\begin{align*}
		&\quad \E\Big[ \sup_{s\in[0,t]} \|  U(s) -U^{(l)}(s) \|^2_{\op} \Big] \nn\\
		&\leq C \Big( \int^t_0\E\big[ \| Df(X(r))  \|^2_{\op} \|  U(r) -U^{(l)}(r) \|^2_{\op} \big]dr + \int^t_0 \E\big[\| Df(X(r))-Df(X^{(l)}(\eta_{N_l}(r)))  \|^2_{\op}\|U^{(l)}(r)\|^2_{\op} \big] dr\Big) \nn \\
		&\leq C  \int^t_0\E\big[ \sup_{s\in[0,r]}  \|  U(s) -U^{(l)}(s) \|^2_{\op} \big]dr +C \frac{1}{N_l}. 
	\end{align*}
	Applying the Gr\"onwall inequality and then taking $t=1$ gives
	\begin{align*} 
		\E\Big[ \sup_{s\in[0,1]} \|  U(s) -U^{(l)}(s) \|^2_{\op} \Big]\leq C \frac{1}{N_l}.
	\end{align*}
	Repeating the same argument for \(U^{-1}\), we also obtain 
	\begin{align}\label{equ_apprUandUlInverse}
		\E\Big[ \sup_{s\in[0,1]} \|  (U(s) )^{-1} -(U^{(l)}(s) )^{-1} \|^2_{\op} \Big]\leq C \frac{1}{N_l}.
	\end{align}
	The uniform boundedness of $G$ and $G^{(l)}$ follows from the uniform boundedness of  $Df$ and $(U^{(l)})^{-1}$, namely,  
	\begin{align}\label{equ_Gbound}
		\sup_{r\in[0,1]}\big( \|G(r)\|_{\op} +\| G^{(l)}(r)\|_{\op} \big)\leq C.
	\end{align}
	Combining \eqref{equ_continuityUlInverse}, \eqref{equ_Dfappro}, \eqref{equ_apprUandUlInverse}, and the uniform boundedness of  $Df$ and $(U^{(l)})^{-1}$,  we have
	\begin{align*}
		&\quad\E\big[ \big\| G(t)- G^{(l)}(t)\big\|_{\op}^2\big]\\
		&\leq C\Big(\E\big[ \big\| \big[(U(t))^{-1} -(U^{(l)}(t))^{-1}\big]Df(X(t)) \big\|^2_{\op}\big] + \E\big[ \big\| \big[(U^{(l)}(t))^{-1}-(U^{(l)}(t_k))^{-1}\big] Df(X(t))\big\|^2_{\op} \big] \\
		&\quad +  \E\big[ \big\|(U^{(l)}(t_k))^{-1} \big(Df(X(t))-Df(X^{(l)}(\eta_{N_l}(t))) \big) \big\|^2_{\op}\big]  \Big)\\
		&\leq C\frac{1}{ N_l}.
	\end{align*} 
	Then for any $\varepsilon>0$,  we have
	\begin{align*}
		\sum^{\infty}_{l=1} 	\PP \Big\{ \int^1_0 \| G(r)- G^{(l)}(r) \|^2_{\op} dr >\varepsilon \Big\} \leq \sum^{\infty}_{l=1} C  \frac{\alpha^{-l}}{\varepsilon} <\infty,
	\end{align*}
	which leads to
	\begin{align}\label{equ_Elimit}
		\lim_{l\to\infty} \int^1_0 \| G(r)- G^{(l)}(r) \|^2_{\op} dr =0,\quad \text{a.s.}
	\end{align}
	
	Based on these estimates, we show the limiting result of $\mathcal V^{(l)}_{M,\sigma}$. 
	Recalling the definition of $\mathcal V_{M,\sigma}$ and  $\mathcal V^{(l)}_{M,\sigma}$ in \eqref{equ_DefVM} and \eqref{equ_DefVMl}, 
	it follows from \eqref{equ_Gbound}, \eqref{equ_Elimit},  and the estimate 
	\begin{align*}
		&\quad \big\| \mathcal V_{M,\sigma}-\mathcal V^{(l)}_{M,\sigma}\big\|_{\op}  \\
		&\leq \frac{1}{12} \int^1_0 \Big[\big\| \big(G(r)- G^{(l)}(r)\big)\sigma \sigma^\top  (G(r))^\top \|_{\op}+\big\| G^{(l)}(r) \sigma \sigma^\top \big((G(r))^\top-(G^{(l)}(r))^\top  \big)\big\|_{\op} \Big]dr\\
		& \leq \frac{\|\sigma\|^2_{\op}}{12} \Big( \int^1_0\big\| G(r)- G^{(l)}(r) \|^2_{\op} dr \Big)^{\frac 12} \Big[ \Big( \int^1_0 \|G(r)\|_{\op}^2 dr\Big)^{\frac 12} +\Big( \int^1_0 \|G^{(l)}(r)\|_{\op}^2 dr\Big)^{\frac 12} \Big]
	\end{align*}
	that
	\begin{align} \label{equ_Vlimit}
		\lim_{l\to\infty} 	\big\| \mathcal V_{M,\sigma}-\mathcal V^{(l)}_{M,\sigma}\big\|_{\op} =0,\quad \text{a.s.}
	\end{align}
	
	Next, for the fixed $v\in\RR^d$, we define 
	\begin{align*}
	\mathbf R^{(l,v)}_{i}(t):=v^\top\big( \mathcal M_{N_i}(t)- \mathcal M^{(l)}_i(t)\big) =\int^t_0 v^\top (G(r)-G^{(l)}(r)) \big(\frac 12-\{N_ir \}\big) \sigma dW(r),\quad t\in[0,1].
	\end{align*}
	Then $\mathbf R^{(l,v)}_i(t)$ is also a martingale and we denote $\mathbf R^{(l,v)}_i:= \mathbf R^{(l,v)}_i(1)$ for simplicity. Therefore, the quadratic variation at $t=1$ satisfies
	\begin{align*}
		\langle \mathbf R^{(l,v)}_{i} \rangle_1 = \int^1_0 \big|   \big((G(r)-G^{(l)}(r) )\sigma\big)^\top v\big|^2 \big(\frac 12-\{N_ir \}\big)^2dr  \leq \frac{\|\sigma\|^2_{\op} |v|^2}{4}  \int^1_0 \| G(r)-G^{(l)}(r)\|^2_{\op}dr.
	\end{align*}
	By defining the event $\mathbf E^{(l,v)}_q:= \{\frac{1}{4}\|\sigma\|^2_{\op} |v|^2  \int^1_0 \| G(r)-G^{(l)}(r)\|^2_{\op}dr \leq q \}$ for $q\in\RR^+$, we have $ 	\langle \mathbf R^{(l,v)}_{i} \rangle_1 \leq q$ on $\mathbf E^{(l,v)}_q$. 
	The Dambis--Dubins--Schwarz time-change theorem for martingales shows that there exists a Brownian motion $ B^{(i,l,v)}$ such that $\mathbf R^{(l,v)}_i = B^{(i,l,v)} ({\langle \mathbf  R^{(l,v)}_i \rangle_1} )$. 
	Therefore, for $\delta>0$,  the reflection principle of Brownian motion gives $\PP\{ \sup_{0\leq t\leq q} |B(t)| >x\} \leq 2\exp\{ -\frac{x^2}{2q}\} $, which leads to
	\begin{align}\label{equ_PofR}
		&\quad \sum^\infty_{i=i_0}\PP\Big\{ | \mathbf R^{(l,v)}_{i}|>(1+\delta)\sqrt{ 2q\log\log N_i},\,\mathbf E^{(l,v)}_q \Big\} \nn \\
		&\leq \sum^\infty_{i=i_0}\PP \Big\{ \sup_{0\leq t\leq q} |B^{(i,l,v)}(t) | >(1+\delta)\sqrt{2q\log \log N_i} \Big\}\nn \\
		& \leq \sum^\infty_{i=i_0}2\exp\big\{ -\frac{(1+\delta)^2 2 q\log\log N_i }{2q}\big\} = \sum^\infty_{i=i_0}2 (i \log \alpha)^{-(1+\delta)^2 }<\infty.
	\end{align}
	Taking $\delta\in\mathbb Q^+$ such that $\delta\to 0$ leads to
	\begin{align}\label{equ_limsuponEq}
		\limsup_{i\to \infty} \frac{| \mathbf R^{(l,v)}_{i}| }{\sqrt{\log \log N_i}} \leq \sqrt{2q} ,\quad \text{ on }  \mathbf E^{(l,v)}_q \backslash F^{(l,v)}_q,
	\end{align}
    where $F^{(l,v)}_q\subset\mathbf E^{(l,v)}_q $ such that $\PP\{F^{(l,v)}_q \}=0$. Let $F^{(l,v)}:=\bigcup_{q\in\mathbb Q^+}F_q^{(l,v)}$. 
	For any $\omega\in\Omega\setminus F^{(l,v)}$, applying \eqref{equ_limsuponEq} with $q\in\mathbb Q^+$ satisfying $q>\frac{\|\sigma\|^2_{\op}|v|^2}{4}  \int^1_0 \| G(r,\omega)-G^{(l)}(r,\omega)\|^2_{\op}dr $ and then letting $q\to  \frac{\|\sigma\|^2_{\op} |v|^2}{4}  \int^1_0 \| G(r,\omega)-G^{(l)}(r,\omega)\|^2_{\op}dr $ yields
	\begin{align*}
		\limsup_{i\to \infty} \frac{| \mathbf R^{(l,v)}_{i}(\omega)| }{\sqrt{\log \log N_i}} \leq  \sqrt{\frac{\|\sigma\|^2_{\op}|v|^2}{2}  \int^1_0 \| G(r,\omega)-G^{(l)}(r,\omega)\|^2_{\op}dr}.
	\end{align*}
    The fact that  $\PP\{\Omega\setminus F^{(l,v)}\}=1$ gives
    \begin{align}\label{equ_limsupR}
		\limsup_{i\to \infty} \frac{| \mathbf R^{(l,v)}_{i}| }{\sqrt{\log \log N_i}} \leq  \sqrt{\frac{\|\sigma\|^2_{\op}|v|^2}{2}  \int^1_0 \| G(r)-G^{(l)}(r)\|^2_{\op}dr},\quad \text {a.s.}
	\end{align}
	Combining \eqref{equ_limsupMl} and \eqref{equ_limsupR}, we obtain
	\begin{align*}
		\limsup_{i\to\infty} \frac{v^\top \mathcal M_{N_i}}{\sqrt{\log\log N_i}}& \geq 	\limsup_{i\to\infty} \frac{v^\top \mathcal M^{(l)}_{i}}{\sqrt{\log\log N_i}} - \limsup_{i \to \infty} \frac{| \mathbf R^{(l,v)}_{i}| }{\sqrt{\log \log N_i}}  \\
		&\geq \sqrt{2v^\top \mathcal V^{(l)}_{M,\sigma} v}  -  \sqrt{\frac{\|\sigma\|^2_{\op}|v|^2}{2}  \int^1_0 \| G(r)-G^{(l)}(r)\|^2_{\op}dr}, \quad \text{a.s.} 
	\end{align*}
	Letting $l\to \infty$, \eqref{equ_Elimit} and \eqref{equ_Vlimit} give
	\begin{align}\label{equ_lowerM}
		\limsup_{i\to\infty} \frac{v^\top \mathcal M_{N_i}}{\sqrt{\log\log N_i}} \geq \sqrt{2v^\top \mathcal V_{M,\sigma} v},\quad \text{a.s.}
	\end{align}
	\subsection{Upper bound estimate in \Cref{prop_LILofM}} \label{subsec_upperLIL}
	The proof of the upper bound does not rely on the construction of a conditionally Gaussian approximation, since the first Borel--Cantelli lemma requires no analysis of dependence within the sequence. In the following, we use martingale estimates to derive the desired result. 
	
	Following an argument similar to \eqref{equ_PofR}, for fixed $v\in\RR^d$, $q>0$, and $\delta>0$, we have 
	\begin{align*}
		&\quad \sum^\infty_{i=1} \PP\Big\{ | v^\top \mathcal M_{N_i}|>(1+\delta)\sqrt{2 q\log\log N_i} ,  \, \langle v^\top \mathcal M_{N_i}\rangle_1 \leq q \Big\} \nn \\
		&\leq \sum^\infty_{i=1}\PP \Big\{ \sup_{0\leq t\leq q} |B^{(i,v)}(t) | >(1+\delta)\sqrt{2 q\log \log N_i} \Big\}\nn \\
		& \leq \sum^\infty_{i=1} 2\exp\big\{ -\frac{(1+\delta)^2  2q\log\log N_i }{2q}\big\} = \sum^\infty_{i=1} 2 (i \log \alpha)^{-(1+\delta)^2 }<\infty.
	\end{align*}
	Therefore, using the Borel--Cantelli lemma gives
	\begin{align}\label{equ_PioforM}
		\PP\Big(\big\{  | v^\top \mathcal M_{N_i}|>(1+\delta) \sqrt{2 q\log\log N_i}\big\} \cap  \big\{\langle v^\top  \mathcal M_{N_i}\rangle_1 \leq q \big\} \quad \text{i.o.} \Big) =0.
	\end{align}
	Here,  we claim that the quadratic variation of the martingale $\{v^\top \mathcal M_{N_i}(t)\}_{t\in[0,1]}$ at $t=1$ 
	\begin{align*}
		\langle  v^\top \mathcal M_{N_i}\rangle_1 = \int^1_0 | (G(r)\sigma)^\top v|^2\big(\frac 12-\{ N_ir\} \big)^2 dr 
	\end{align*}
	satisfies 
	\begin{align}\label{equ_limitvariationVM}
		\lim_{i\to \infty}	\langle v^\top  \mathcal M_{N_i}\rangle_1   =v^\top \mathcal V_{M,\sigma} v, \quad \text{a.s.}
	\end{align}
	Next we prove \eqref{equ_limitvariationVM}.
	Note that
	\begin{align}\label{equ_hNL2}
		&\quad \lim_{i\to \infty}  \int^1_0 \big|(G(r)\sigma)^\top v\big|^2 \big(\frac 12-\{N_i r\} \big)^2   dr \nn\\
		&=	\lim_{i\to \infty}\sum^{N_i-1}_{k=0}  \int^{\frac{k+1}{N_i}}_{\frac{k}{N_i}} \big|(G(r)\sigma)^\top v \big|^2 \big(\frac 12-(N_ir-k) \big)^2 dr  \nn\\
		&=\lim_{i\to \infty}\sum^{N_i-1}_{k=0} \Big[ \big|( G(\frac{k}{N_i})\sigma)^\top v \big|^2 \int^{\frac{k+1}{N_i}}_{\frac{k}{N_i}}   \big(\frac 12-(N_ir-k) \big)^2 dr \nn \\
		&\quad  +\int^{\frac{k+1}{N_i}}_{\frac{k}{N_i}}   v^\top \big[ G(r)\sigma\sigma^\top(G(r))^\top - G(\frac{k}{N_i})\sigma\sigma^\top(G(\frac{k}{N_i}))^\top\big] v \big(\frac 12-(N_ir-k) \big)^2 dr\Big]  \nn\\
		&=\frac{1}{12} \lim_{i \to \infty}\frac{1}{N_i}\sum^{N_i-1}_{k=0}\big|( G(\frac{k}{N_i})\sigma)^\top v \big|^2 \nn \\
		&\quad + \lim_{i\to \infty}\sum^{N_i-1}_{k=0} \int^{\frac{k+1}{N_i}}_{\frac{k}{N_i}}  v^\top  \big[ G(r)\sigma\sigma^\top(G(r))^\top - G(\frac{k}{N_i})\sigma\sigma^\top(G(\frac{k}{N_i}))^\top\big] v \big(\frac 12-(N_ir-k) \big)^2 dr.
	\end{align}
	It follows from \eqref{equ_XMomentcontinuity} 
	that $G$ has continuous sample paths almost surely, which shows that the corresponding Riemann sums converge almost surely.  Therefore, we arrive at
	\begin{align*}
		\frac{1}{12} \lim_{i\to \infty} \frac{1}{N_i} \sum^{N_i-1}_{k=0} G(\frac{k}{N_i})\sigma\sigma^\top(G(\frac{k}{N_i}))^\top = 	\frac{1}{12} \int^1_0 G(r)\sigma\sigma^\top(G(r))^\top dr,\quad \text{a.s.}
	\end{align*}
	and
	\begin{align*}
		\lim_{i\to \infty}\sum^{N_i-1}_{k=0} \int^{\frac{k+1}{N_i}}_{\frac{k}{N_i}}    \big[ G(r)\sigma\sigma^\top(G(r))^\top - G(\frac{k}{N_i})\sigma\sigma^\top(G(\frac{k}{N_i}))^\top\big]   \big(\frac 12-(N_ir-k) \big)^2 dr=\textbf 0_{d\times d}, \quad \text{a.s.}
	\end{align*}
	Plugging these two estimates into \eqref{equ_hNL2}  proves \eqref{equ_limitvariationVM}.

Fix $v\in\RR^d$, $q>0$.
    It follows from \eqref{equ_PioforM} and \eqref{equ_limitvariationVM} that  for $\delta >0$, there exists a null set $F_{\delta}\subset \{ v^\top \mathcal V_{M,\sigma}(\omega) v< q\}$ such that for every $\omega\in \{ v^\top \mathcal V_{M,\sigma}(\omega) v< q\}\backslash F_{\delta}$,  $\langle v^\top  \mathcal M_{N_i}\rangle_1\leq q$ and 
    $| v^\top \mathcal M_{N_i}|\leq (1+\delta) \sqrt{2 q\log\log N_i} $
    for any sufficiently large $i\in\N^+$.
    This yields
	\begin{align*}
		\limsup_{i \to \infty} \frac{\big| v^\top \mathcal M_{N_i}\big| }{ \sqrt{\log\log N_i}} \leq (1+\delta) \sqrt{2q},\quad \text{a.s. on } \{v^\top \mathcal V_{M,\sigma}  v<q\},
	\end{align*}
    which further gives
    	\begin{align*}
		\limsup_{i \to \infty} \frac{\big| v^\top \mathcal M_{N_i}\big| }{ \sqrt{\log\log N_i}} \leq  \sqrt{2q},\quad \text{a.s. on } \{v^\top \mathcal V_{M,\sigma}  v<q\},
	\end{align*}
    by taking $\delta\in\mathbb Q^+$ such that $\delta\to 0$. 
	Hence, by an argument analogous to that leading to \eqref{equ_limsupR}, we obtain
	\begin{align} \label{equ_upperM}
		\limsup_{i \to \infty} \frac{\big| v^\top \mathcal M_{N_i}\big|}{\sqrt{\log\log N_i}} \leq \sqrt{2 v^\top \mathcal V_{M,\sigma} v},\quad\text{a.s.},
	\end{align}
	
	Combining \eqref{equ_lowerM} and \eqref{equ_upperM}, we complete the proof of \Cref{prop_LILofM}.
	
	\section{LIL constants and numerical comparison of methods}\label{sec_Comparison}

In this section, we compare numerical methods through their LIL constants and illustrate the theoretical results numerically. We first compare the \(\theta\)-Euler--Maruyama methods, which share the same LIL constants, and examine their differences at finite refinement levels. We then compare these methods with the BAOAB method to illustrate how a different discretization of the stochastic term can lead to different LIL constants.

	\subsection{Common LIL constants and finite-level differences
    }\label{subsec_thetaEM}
	
We use a linear scalar equation to illustrate how the \(\theta\)-Euler--Maruyama methods can share the same LIL constants while exhibiting different error behavior at finite refinement levels. Specifically, we consider
		\begin{align} \label{equ_OU}
		dX(t) = \mu X(t) dt +\sigma dW(t), \quad X(0)=x_0,
	\end{align} 
    where $t\in[0,1]$, $\mu\neq 0$, and $\sigma>0$.
Let \(X^{N,\theta}\) denote its numerical approximation obtained by the \(\theta\)-Euler–Maruyama method with step size \(N^{-1}\). Specializing  \Cref{thm_thetamethod} to this equation gives  
		\begin{align*}
			\limsup_{i\to \infty} \frac{N_i(X(1)-X^{N_i,\theta}(1))}{\sqrt{\log\log N_i} }&=\sigma \sqrt{\frac{1}{12} \mu\big(\exp(2\mu)-1\big) },\quad \text{a.s.}\\
			\liminf_{i\to \infty} \frac{N_i(X(1)-X^{N_i,\theta}(1))}{\sqrt{\log\log N_i} }&=-\sigma \sqrt{\frac{1}{12} \mu\big(\exp(2\mu)-1\big) },\quad \text{a.s.}
		\end{align*}
        Although the LIL constants are independent of \(\theta\), the choice of \(\theta\) can substantially affect the normalized errors at finite refinement levels.

		To explain this point, we examine the asymptotic behavior of the difference between two numerical solutions. For any fixed $\theta_1,\theta_2\in[0,1]$, we have
		\begin{align}\label{equ_asO1}
			\lim_{i\to\infty}N_i
			\bigl(X^{N_i,\theta_1}(1)-X^{N_i,\theta_2}(1)\bigr)
			&=(\theta_1-\theta_2)
			\Bigl[
			\mu^2e^\mu x_0
			+\mu\sigma\int_0^1e^{\mu(1-s)}
			\bigl[1+\mu(1-s)\bigr]\,dW(s)
			\Bigr]\nonumber\\
			&=:C^{\theta_1,\theta_2},
			\quad\text{a.s.}
		\end{align}
		See \Cref{subsec_proofinseccompare} for the proof of \eqref{equ_asO1}.
        Define
\[
\mathcal E^{\theta_1,\theta_2}_{N_i}
:=
\frac{N_i\left|X^{N_i,\theta_1}(1)-X^{N_i,\theta_2}(1)\right|}
{\sqrt{\log\log N_i}}.
\] 
Recall that $N_i=\alpha^i$, $i\in\N^+$. It follows from \eqref{equ_asO1} that for sufficiently large $i\in\N^+$,
		$$
		\mathcal E^{\theta_1,\theta_2}_{N_i}
		=\frac{|C^{\theta_1,\theta_2}|+o(1)}
		{\sqrt{\log(i\log \alpha)}},
		\quad\text{a.s.}
		$$
Consequently, these normalized differences converge to zero very slowly whenever \(C^{\theta_1,\theta_2}\neq0\).

			\begin{figure}
			\begin{center}
				\includegraphics[scale=0.20]{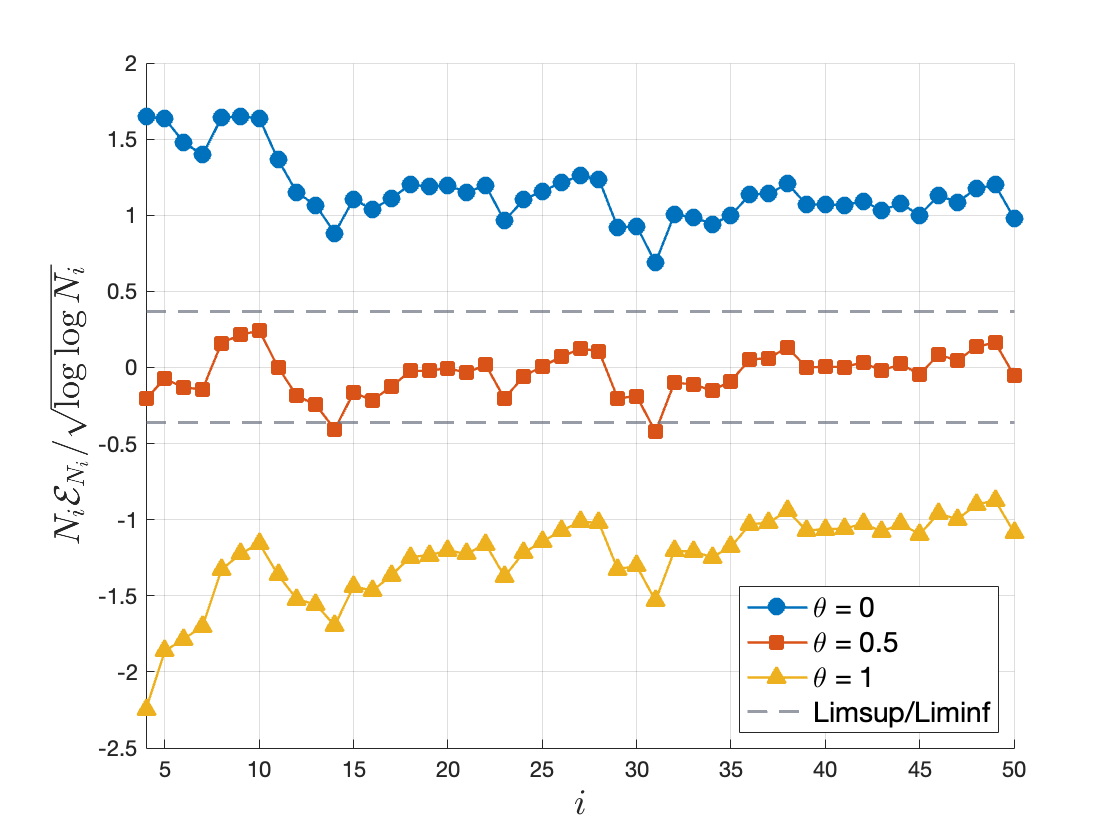}
			\end{center}
			\caption{Comparison of $\theta$-methods}
			\label{fig:theta_comparison}
		\end{figure}
        
 We illustrate this behavior with $\mu=1$, $\sigma=0.5$,  $X(0)=1$, and $N_i=2^i$ for $i=4,\ldots, 50$. The corresponding LIL constants are $ \pm\sigma \sqrt{\frac{1}{12} \mu\big(\exp(2\mu)-1\big) } \approx \pm 0.3648$. 
 Figure~\ref{fig:theta_comparison} displays the normalized errors 
		$$
		\frac{N_i(X(1)-X^{N_i,\theta}(1))}{\sqrt{\log\log N_i}}
		$$
		for $\theta=0,\frac12,1$, computed along the same sample path. 
		The dashed lines indicate the theoretical LIL bounds. Over the displayed range, the midpoint method, corresponding to \(\theta=\frac12\), exhibits fluctuations around zero, whereas the other two methods show substantial offsets. These offsets are consistent with the slowly vanishing differences described above. 
\begin{table}[htbp]
    \centering
    \caption{Normalized differences between the $\theta$-Euler--Maruyama methods}
    \label{tab:theta_difference}
    \renewcommand{\arraystretch}{1.15}
    \setlength{\tabcolsep}{8pt}
    \begin{tabular}{@{}lcccccc@{}}
        \toprule
         
        & $N_i=2^{5}$ & $N_i=2^{10}$ & $N_i=2^{20}$
        & $N_i=2^{30}$ & $N_i=2^{40}$ & $N_i=2^{50}$ \\
        \midrule
        $\mathcal E^{1,1/2}_{N_i}$
        & 1.7903 & 1.3980 & 1.1990 & 1.1160 & 1.0666 & 1.0325 \\
        $\mathcal E^{0,1/2}_{N_i}$
        & 1.7068 & 1.3959 & 1.1990 & 1.1160 & 1.0666 & 1.0325 \\
        \bottomrule
    \end{tabular}
\end{table}

 Table~\ref{tab:theta_difference} reports the normalized differences between the explicit and implicit Euler–Maruyama methods and the midpoint method. The differences decrease slowly and remain substantial even at the finest refinement level shown. Indeed, a leading-order estimate for the refinement index required to reduce \(\mathcal E^{\theta_1,\theta_2}_{N_i}\) below \(0.1\) is
$
i\gtrsim
\frac{\exp(100|C^{\theta_1,\theta_2}|^2)}
{\log2}.
$ 
For the sample path considered here,
$|C^{1,1/2}|=|C^{0,1/2}|\approx1.9441,$
which gives an estimated refinement index of order \(10^{164}\). This illustrates why convergence of the normalized differences to zero is difficult to observe over a computationally accessible range.

\begin{figure}
			\begin{center}
				\includegraphics[scale=0.20]{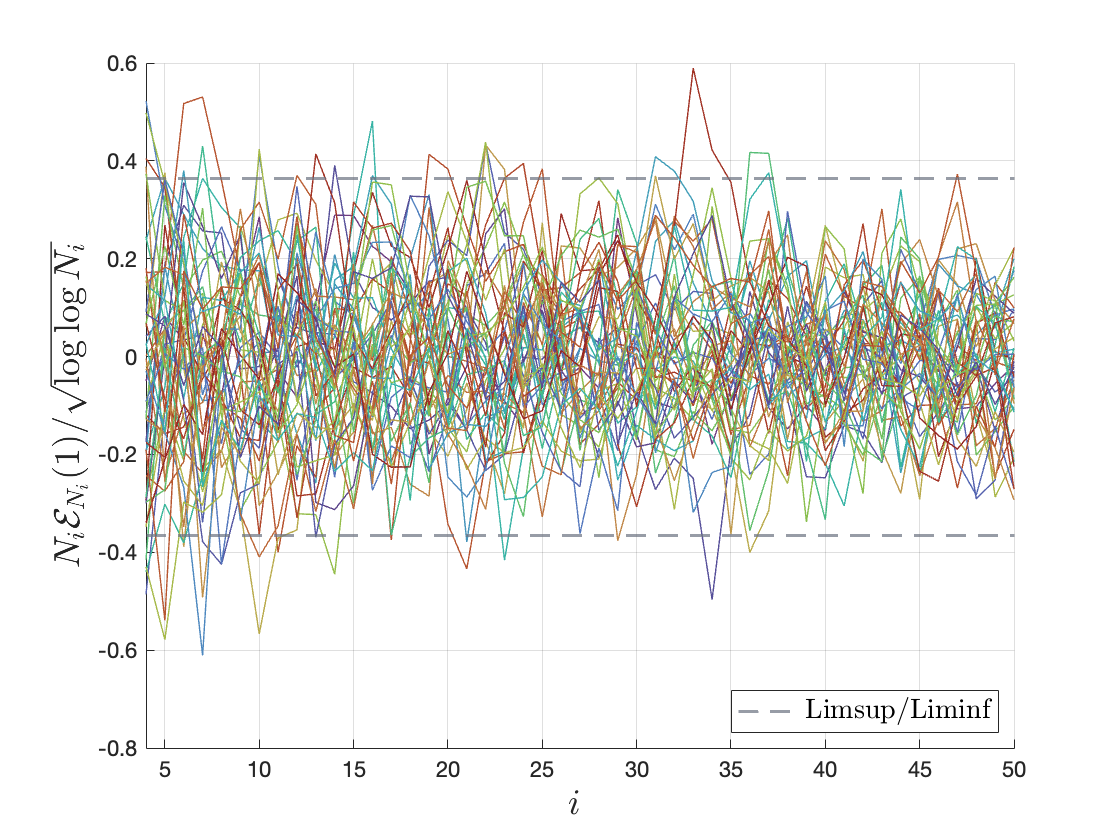}
			\end{center}
			\caption{Almost sure asymptotics of the midpoint method ($\theta=1/2$)}
			\label{fig:midpoint}
		\end{figure}

To further illustrate that the observed behavior of the midpoint method is not specific
to a single sample path, Figure \(\ref{fig:midpoint}\) displays the normalized errors of the midpoint method along \(50\) sample paths. The observed fluctuations are consistent with the theoretical LIL bounds, supporting the predicted almost sure asymptotic behavior. 

		
        
		\subsection{Method-dependent LIL constants} \label{subsec_BAOAB}
		The common LIL constants of the \(\theta\)-Euler--Maruyama methods reflect the structure of this family of discretizations. Changing the
		discretization of the stochastic term, however, may lead to different LIL
		constants. 
To illustrate this point, we compare the almost sure error
asymptotics of the $\theta$-Euler--Maruyama methods and the
BAOAB method (see, e.g., \cite[Section 7.3]{Leimkuhler2015book})
for the linear underdamped Langevin equation on $[0,1]$
\begin{equation}\label{eq:origin}
		\left\{\begin{aligned}
			dP(t) &= -\mu P(t) dt-\gamma Q(t)  dt+\lambda dW(t) , \quad P(0)=p_0,
			\\
			dQ(t)  &= P(t) dt, \quad Q(0)=q_0,
		\end{aligned} \right. 
	\end{equation}
	where  $\mu,\gamma, \lambda >0$ are constants.  Setting 
	\begin{align*}
		X(t)= \begin{pmatrix}
			P(t) \\ Q(t)
		\end{pmatrix}, \quad A= \begin{pmatrix}
			-\mu & -\gamma \\ 1 & 0
		\end{pmatrix}, \quad \sigma = \begin{pmatrix}
			\lambda \\ 0
		\end{pmatrix},
	\end{align*}
and $f(X(t))=AX(t)$,	we rewrite \eqref{eq:origin} in the form of  \eqref{equ_nonliearGeneral}. 
    Then the LILs for these methods take the same form as those in \Cref{thm_LILadditive},  with 
\begin{align*}
			\mathcal V_\sigma 
			&= \frac{\lambda^2}{12}  \int^1_0 e^{At}  \begin{pmatrix}
				\mu^2 & -\mu \\ -\mu & 1
			\end{pmatrix}   (e^{At})^\top  dt
		\end{align*}
for the $\theta$-Euler--Maruyama method, and
	\begin{align*}
		\mathcal V^S_\sigma 
		&= \frac{\lambda^2}{12}  \int^1_0 e^{At}  \begin{pmatrix}
			0 & 0 \\ 0 & 1
		\end{pmatrix}   (e^{At})^\top  dt
	\end{align*}
for the BAOAB method.
    The proofs are provided in Section~\ref{subsec_proofinseccompare}.  These results also demonstrate the applicability of our framework to other numerical methods, with the leading error term and the corresponding LIL constants determined by the particular discretization.

We take 
\(\mu=2\) and \(\gamma=1\) as an example to evaluate the LIL constants explicitly, which are provided in Table~\ref{tab:LILconstants}. Here, $\mathcal E^P_{N_i}(1)$ and $\mathcal E^Q_{N_i}(1)$
denote the errors in the $P$- and $Q$-components, respectively,
between the exact solution and the corresponding numerical
solution, and $|\mathcal E_{N_i}(1)|$ denotes the Euclidean
norm of the error. The constants in the table are the almost
sure limsup values of these quantities after multiplication
by $N_i/\sqrt{\log\log N_i}$. Although the two methods have the same strong convergence
order, their LIL constants differ, revealing distinct
almost sure error envelopes under mesh refinement. 

\begin{table}[htbp] 
    \centering
    \caption{ LIL constants for $\mu=2$ and $\gamma=1$}
    \label{tab:LILconstants}
    \renewcommand{\arraystretch}{1.3}
    \setlength{\tabcolsep}{8pt}
    \begin{tabular}{@{}lcc@{}}
        \toprule
        Error quantity
        & $\theta$-Euler--Maruyama
        & BAOAB  \\
        \midrule
        $\mathcal E^P_{N_i}(1)$
        & $
           \frac{\lambda}{2\sqrt{6}}\sqrt{5-e^{-2}}$
        & $\frac{\lambda}{2\sqrt{6}}\sqrt{1-5e^{-2}}$
        \\[6pt]
        $\mathcal E^Q_{N_i}(1)$
        & $\frac{\lambda}{2\sqrt{6}}\sqrt{1-e^{-2}}$
        & $\frac{\lambda}{2\sqrt{6}}\sqrt{5-13e^{-2}}$
        \\[6pt]
        $|\mathcal E_{N_i}(1)|$
        & $\frac{\lambda}{2\sqrt{6}}
           \sqrt{3-e^{-2}+2\sqrt{2}}$
        & $\frac{\lambda}{2\sqrt{6}}
           \sqrt{3-9e^{-2}+2\sqrt{2-12e^{-2}+20e^{-4}}}$
        \\[6pt]
        \bottomrule
\end{tabular}
\end{table}

		
We now numerically illustrate these almost sure error
asymptotics. Keeping $\mu=2$ and $\gamma=1$, we set $\lambda=1$ and
$(P(0),Q(0))^\top=(1,1)^\top$,  and take
$N_i=2^i$ for $i=4,\ldots,30$.
We use the midpoint method as a representative of the
$\theta$-Euler--Maruyama family.
 
	\begin{figure}[ht]
		\begin{subfigure}{0.3\textwidth}
			\centerline{
				\includegraphics[scale=0.17]{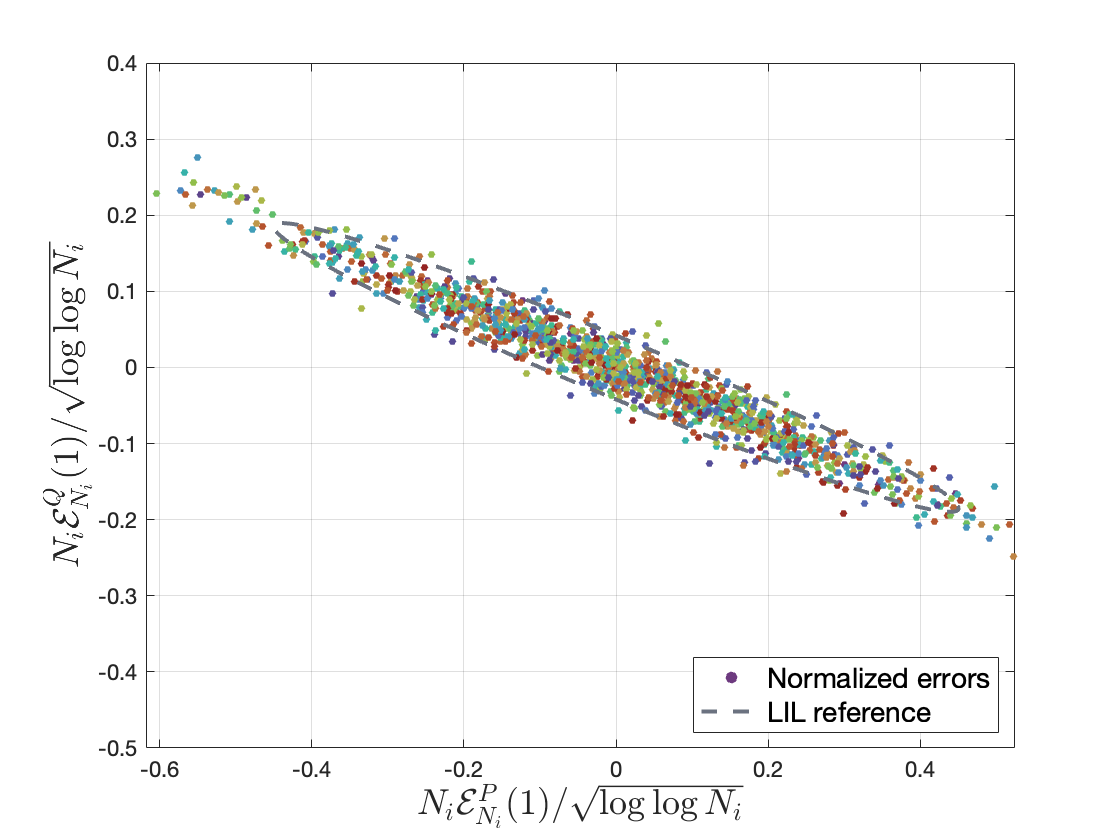}}
			\caption{Midpoint method}
		\end{subfigure}
		\hspace{30mm}
		\begin{subfigure}{0.3\textwidth}
			\centerline{
				\includegraphics[scale=0.17]{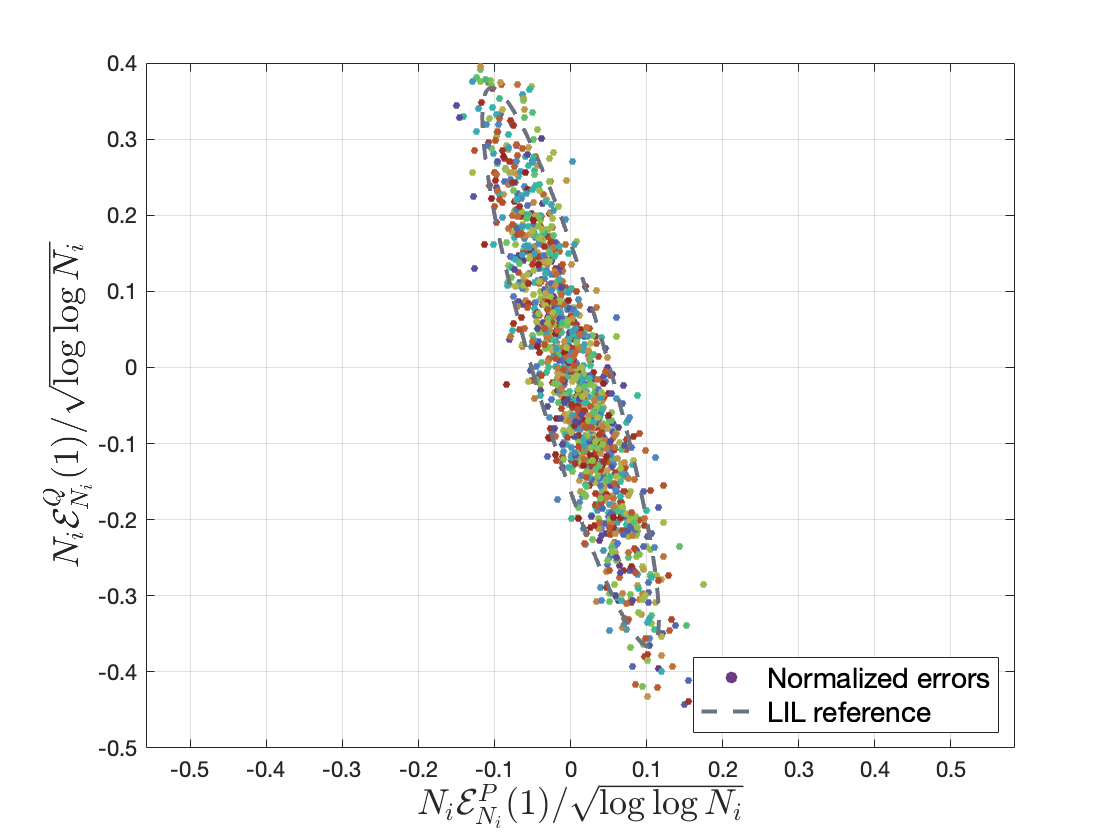}}
			\caption{BAOAB method}
		\end{subfigure}
		\caption{Directional almost sure asymptotics}
		\label{fig:langevinComparison2d}
	\end{figure}
	
	We first illustrate the directional behavior of the errors. In 
	Figure~\ref{fig:langevinComparison2d}(A) and 	Figure~\ref{fig:langevinComparison2d}(B), we display the normalized error vectors $
	\big(\frac{N_i\mathcal E^P_{N_i}(1)}
	{\sqrt{\log\log N_i}},\frac{N_i\mathcal E^Q_{N_i}(1)}
	{\sqrt{\log\log N_i}} \big)
	$
	along $50$ sample paths for the midpoint and BAOAB methods, respectively. 
    Each point 
    corresponds to a refinement index $i$, and
	the dashed reference curves are determined by the directional LIL constants. The plotted points lie almost entirely within the region determined by the dashed curves, and some reach its boundary. These observations agree with the theoretical results and numerically demonstrate the sharpness of the bounds.
	
	\begin{figure}[ht]
		\begin{subfigure}{0.3\textwidth}
			\centerline{
				\includegraphics[scale=0.17]{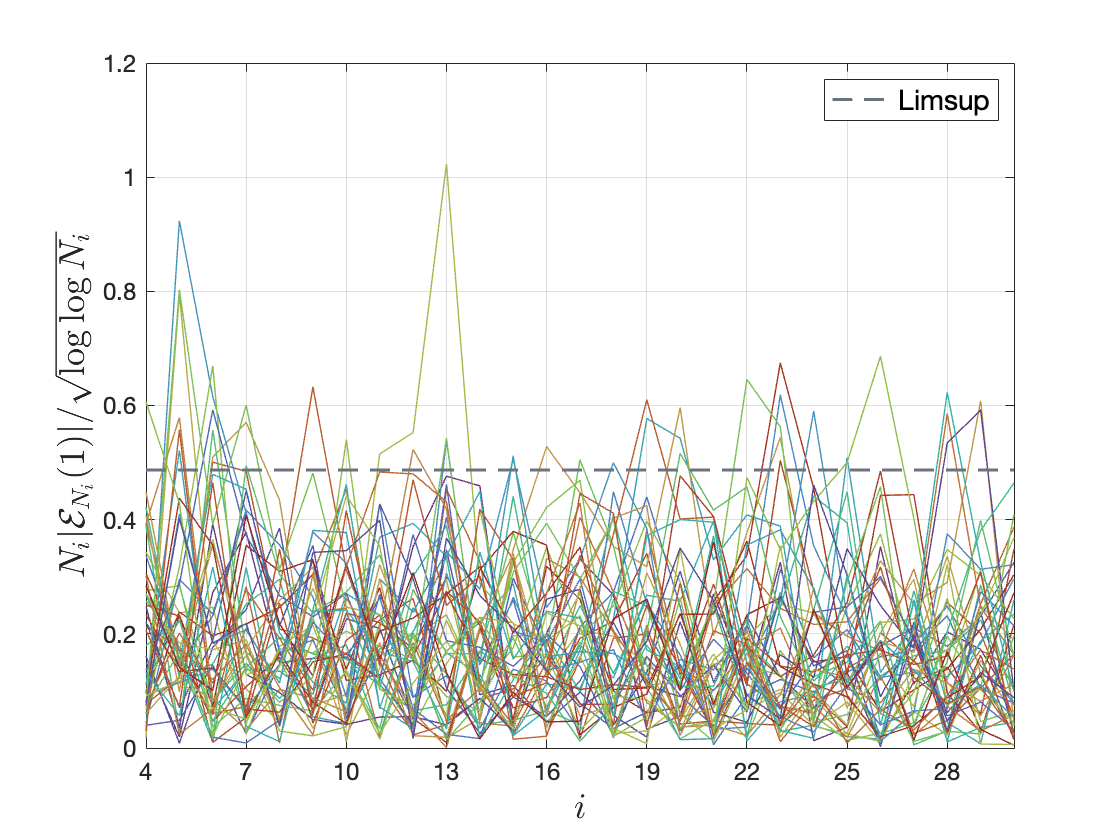}}
			\caption{Midpoint method}
		\end{subfigure}
		\hspace{30mm}
		\begin{subfigure}{0.3\textwidth}
			\centerline{
				\includegraphics[scale=0.17]{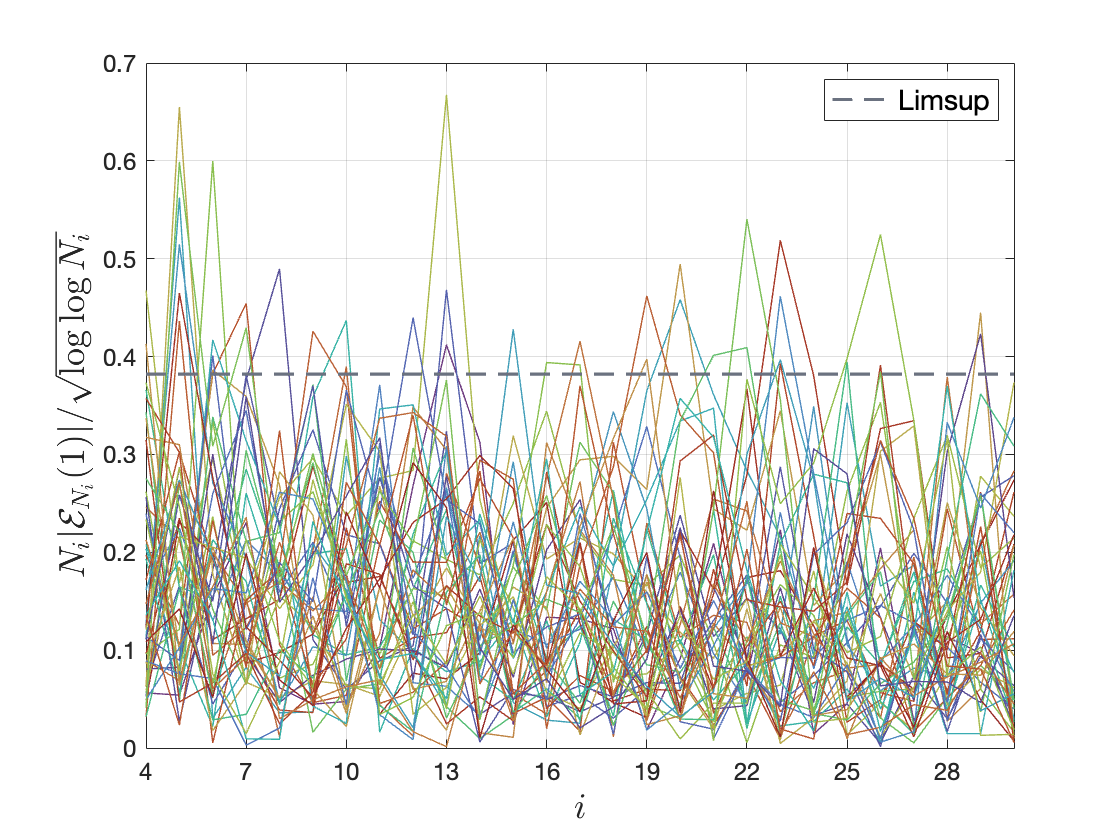}}
			\caption{BAOAB method}
		\end{subfigure}
		\caption{Almost sure asymptotics of the norm}
		\label{fig:langevinComparison}
	\end{figure}
	We next illustrate the almost sure asymptotic behavior of the error norm. Figure~\ref{fig:langevinComparison}(A)--(B) display the 
	norms of the normalized errors along 50 sample paths for the midpoint and BAOAB methods, respectively. 
	The dashed lines indicate the theoretical LIL constants,  approximately $0.4870$ for the $\theta$-methods and $0.3822$ for the BAOAB method, respectively. The observed fluctuations are also consistent with the theoretical LIL bounds. Most of the normalized error norms lie below the theoretical limsup levels, and some reach these levels, providing numerical evidence of the sharpness of the bounds.

	
	\appendix
	\section{Supplement}
			\subsection{Proof of \Cref{thm_thetamethod} }\label{subsec_prooftheta}
	Fix $\theta\in[0,1]$. Throughout this proof, we retain the notation $X_k$ for the numerical solution, $X^N(t)$ for its continuous-time extension, and $\mathcal E_N(t):=X(t)-X^N(t)$ for the error process, omitting the parameter $\theta$ for simplicity.
	The continuous-time extension is given by
	$$X^N(t) = X_0 +(1-\theta) \int_0^t f(X^N(\eta_N(s)))ds+\theta  \int_0^t f(X^N(\zeta_N(s)))ds  + \sigma W(t).$$
	Subtracting this expression from the integral equation for \(X(t)\) and applying the Taylor expansion, we obtain
    \begin{align}\label{equ_thetaENTaylor}
		\mathcal{E}_N(t) &=(1-\theta) \int_0^t [f(X(s)) - f(X^N(\eta_N(s)))] ds+\theta \int_0^t [f(X(s)) - f(X^N(\zeta_N(s)))] ds \nn \\
		&=(1-\theta) \int_0^t\Big(  Df(X(s))(X(s) - X^N(\eta_N(s))) - \bar{\xi}_{1,N}(s) (X(s) - X^N(\eta_N(s)), \nn\\
		&\qquad X(s) - X^N(\eta_N(s)) ) \Big)ds \nn\\
		&\quad + \theta  \int_0^t\Big(  Df(X(s))(X(s) - X^N(\zeta_N(s))) - \bar{\xi}_{2,N}(s) (X(s) - X^N(\zeta_N(s)),  \nn\\
		&\qquad X(s) - X^N(\zeta_N(s)) )  \Big)ds  \nn\\
		&=  \int_0^t Df (X(s))\big[(1-\theta)(X(s) - X^N(\eta_N(s)) ) +\theta(X(s) - X^N(\zeta_N(s)) ) \big]  ds\nn \\
		&\quad- (1-\theta) \int_0^t \bar{\xi}_{1,N}(s) (X(s) - X^N(\eta_N(s)),X(s) - X^N(\eta_N(s)) ) ds \nn \\
		&\quad - \theta  \int_0^t \bar{\xi}_{2,N}(s) (X(s) - X^N(\zeta_N(s)),X(s) - X^N(\zeta_N(s)) ) ds ,
	\end{align}
	where 	$\bar{\xi}_{1,N}(s) (u,v) :=\int_0^1 \varepsilon_{1} 
	D^2f\big( \varepsilon_{1} 
	X(s)+ (1-\varepsilon_{1} ) X^N(\eta_N(s))
	\big) (u,v) d\varepsilon_{1} $ and $\bar{\xi}_{2,N}(s) (u,v) :=\int_0^1 \varepsilon_{2} 
	D^2f\big( \varepsilon_{2} 
	X(s)+ (1-\varepsilon_{2} ) X^N(\zeta_N(s))
	\big) (u,v) d\varepsilon_{2} $
	for $u,v\in\RR^d$.
    
	Using the identities
	\begin{align*}
		&\quad X^N(s) - X^N(\eta_N(s)) \\
		&=(1-\theta)   f(X^N(\eta_N(s)))(s - \eta_N(s)) + \theta   f(X^N(\zeta_N(s)))(s - \eta_N(s)) + \sigma\big(W(s)- W(\eta_N(s))\big)
	\end{align*}
	and
	\begin{align*}
		&\quad X^N(s) - X^N(\zeta_N(s)) \\
		&=(1-\theta)   f(X^N(\eta_N(s)))(s - \zeta_N(s)) + \theta   f(X^N(\zeta_N(s)))(s - \zeta_N(s))+ \sigma\big(W(s)- W(\zeta_N(s))\big),
	\end{align*}
	we obtain
	\begin{align*}
		&\quad (1-\theta) \big[	X(s) - X^N(\eta_N(s))\big] +\theta  \big[	X(s) - X^N(\zeta_N(s))\big]\\
		&= 	\mathcal{E}_N(s)+ 
		(1-\theta) \big[	X^N(s) - X^N(\eta_N(s))\big] +\theta  \big[	X^N(s) - X^N(\zeta_N(s))\big]  \\
		&=	\mathcal{E}_N(s)+ (1-\theta)\Big[(1-\theta)   f(X^N(\eta_N(s)))(s - \eta_N(s)) + \theta   f(X^N(\zeta_N(s)))(s - \eta_N(s)) \\
		&\quad + \sigma\big(W(s)- W(\eta_N(s))\big) \Big] + \theta \Big[ (1-\theta)   f(X^N(\eta_N(s)))(s - \zeta_N(s)) \\
		&\quad + \theta   f(X^N(\zeta_N(s)))(s - \zeta_N(s))+ \sigma\big(W(s)- W(\zeta_N(s))\big)\Big]\\
		&= 	\mathcal{E}_N(s)+ \Big[(1-\theta)   f(X^N(\eta_N(s)))+ \theta   f(X^N(\zeta_N(s))) \Big]\Big[(1-\theta) (s - \eta_N(s)) +\theta(s-\zeta_N(s)) \Big]\\
		&\quad + \sigma \big[ (1-\theta)(W(s)-W(\eta_N(s)))+\theta (W(s)- W(\zeta_N(s))) \big] .
	\end{align*}
	Plugging this into the first term on the right-hand side of \eqref{equ_thetaENTaylor}
	and using $\zeta_N(s)=\eta_N(s)+\frac 1N$, we have
	\begin{align*}
		&\quad  \int_0^t Df(X(s))\big[(1-\theta)(X(s) - X^N(\eta_N(s)) ) +\theta(X(s) - X^N(\zeta_N(s)) ) \big]  ds\\
		&=\int_0^t Df(X(s)) 	\mathcal{E}_N(s)  ds \\
		&\quad + \int_0^t Df(X(s))\big[ (1-\theta)f(X^N(\eta_N(s)))+\theta f(X^N(\zeta_N(s))) \big] \Big[(s - \eta_N(s)) -\frac \theta N\Big] ds \\
		&\quad+ \int_0^t  Df(X(s))\sigma \big( (1-\theta)(W(s)-W(\eta_N(s)))+ \theta(W(s)- W(\zeta_N(s))) \big) ds.
	\end{align*}
	An argument similar to that in \Cref{sec_preEstimates} shows that 
	the dominant contributions to the terminal error come from the following two integrals:
    $$(1-\theta) \int^1_0 \mathcal J(1,s)Df(X(s))  \sigma(W(s) - W(\eta_N(s))) ds$$ and $$\theta \int^1_0\mathcal J(1,s)Df(X(s))  \sigma(W(s) - W(\zeta_N(s))) ds.$$
The first term has already been treated in \(\eqref{equ_FubiniMain}\). For the second term, the product formula for the divergence operator and the Fubini theorem yield  	\begin{align}\label{equ_FUbiniSecond}
		&\quad	\int^1_0  \mathcal J(1,s)Df(X(s)) \sigma(W(s) - W(\zeta_N(s))) ds \nn\\
		&= \int^1_0g(s)\Big( \int^s_{\zeta_N(s)} \sigma  dW(r)\Big) ds\nn \\
		&=  \int^1_0  \int^{\eta_N(r)}_{r} g(s)\sigma ds \delta W(r)  + \sum^m_{j=1} \int^1_0 \int^s_{\zeta_N(s)} (\mathcal D_r g(s)[e_j]) (\sigma e_j) dr ds\nn\\
		&=  \int^1_0  g(r) (\eta_N(r)-r)  \sigma \delta W(r) + \int^1_0  \int^{\eta_N(r)}_{r}  (g(s)-g(r))  \sigma ds\delta W(r) \nn  \\
		&\quad + \sum^m_{j=1} \int^1_0 \int^s_{\zeta_N(s)} (\mathcal D_r g(s)[e_j]) (\sigma e_j) dr ds.
	\end{align}
	Combining \eqref{equ_FubiniMain} and \eqref{equ_FUbiniSecond}, we obtain the following decomposition of the leading Skorohod integral in the rescaled terminal error \(N\mathcal E_N(1)\) for the $\theta$-Euler--Maruyama method
	\begin{align*}
		&\quad	N \int^1_0  \mathcal J(1,r)Df(X(r)) \big[ (1-\theta )(\zeta_N(r)-r) +	\theta	 (\eta_N(r)-r) \big]   \sigma \delta W(r) \\
		&=N  \int^1_0  \mathcal J(1,r)Df(X(r))  (\zeta_N(r)-r)\sigma  \delta W(r) - \theta \int^1_0  \mathcal J(1,r)Df(X(r))  \sigma  \delta W(r)\\
        &=Y_{N}+\big(\frac12-\theta\big) \int_0^1 \mathcal J(1,r)Df(X(r))\sigma  \delta W(r),
	\end{align*}
where $Y_N$ is	defined in \eqref{equ_YNexpress}. The second term is independent of \(N\) and finite almost surely, so it vanishes after division by \(\sqrt{\log\log N_i}\) as \(i\to\infty\). Consequently, the LIL asymptotics are still determined by \(Y_{N_i}\), and the resulting LIL constants coincide with those in \Cref{thm_LILadditive}.

    \subsection{Proof of Proposition~\ref{prop_gDrg}} \label{subsec_proofofprop_gDrg}
              Fix $p\geq 2$.     Without loss of generality, assume that $0\leq t\leq s\leq 1$ in the proof. 

                   We first prove \eqref{equ_Lpofg}. 
		The Lipschitz continuity of $f$ and the uniform  moment boundedness of $X$ yield
		\begin{align}\label{equ_XMomentcontinuity}
			\big(	\E[| X(s)-X(t)|^{p} ] \big)^{\frac 1p} &\leq  \Big( \E \Big[\big| \int^s_t f(X(r))dr \big|^{p} \Big] \Big)^{\frac 1p}+ \|\sigma\|_{\op} \Big(\E \Big[\big|\int^s_t dW(r)\big|^p \Big] \Big)^{\frac 1p}  \leq  C_p  (s-t)^{\frac 12}.
		\end{align}
		Combining
		\begin{align*}
			\Big(	\E\Big[\big\| g(s)-g(t)\big\|_{\op}^{p} \Big] \Big)^{\frac 1p} 
			&\leq  \Big(	\E\Big[\big\| \big(\mathcal J(1,s) -\mathcal J(1,t)  \big) Df(X(s)) \big\|_{\op}^{p}\Big] \Big)^{\frac 1p}\\
			&\quad + \Big(	\E\Big[\big\|\mathcal J(1,t) \big(  Df(X(s))- Df(X(t))\big) \big\|_{\op}^{p} \Big] \Big)^{\frac 1p}
		\end{align*}
		with  \eqref{equ_JTemporalHolder}, \eqref{equ_XMomentcontinuity},  and the uniform boundedness of  $\|Df(\cdot)\|_{\op}$, we have
		\begin{align*} 
			\Big(	\E\Big[\big\| g(s)-g(t)\big\|_{\op}^{p} \Big] \Big)^{\frac 1p} \leq C_p |s-t|^{\frac 12}.
		\end{align*}

        Next we prove \eqref{equ_LpofDg}.
		Taking the Malliavin derivative of \eqref{equ_nonliearGeneral} gives 
		\begin{align*}
			\mathcal D_r X(z) = \int_0^z Df(X(t)) \mathcal D_r X(t) dt + \mathbf{1}_{\{r \leq z\}} \sigma,
		\end{align*}
		which leads to 
		\begin{align} \label{equ_DrX}
			\mathcal D_r X(z) = \mathbf{1}_{\{r \leq z\}}\mathcal J(z,r)\sigma  =\mathbf{1}_{\{r \leq z\}} U(z)(U(r))^{-1}  \sigma.
		\end{align}
Using $\mathcal J(s,t)=U(s)(U(t))^{-1}$,  \eqref{equ_DrU},   and
$$
\mathcal D_r((U(t))^{-1})[u]
=-(U(t))^{-1}(\mathcal D_rU(t)[u])(U(t))^{-1},
$$
for $r\leq s$ we obtain
\begin{align*}
(\mathcal D_r\mathcal J(s,t)[u])v
&=(\mathcal D_rU(s)[u])(U(t))^{-1}v
-\mathcal J(s,t)(\mathcal D_rU(t)[u])(U(t))^{-1}v\\
   &= \int^s_t \mathcal J(s,z) \big(\mathcal D_r \big(Df(X(z))\big) [u] \big)(\mathcal J(z,t)v) dz  \\
			&= \int^s_t \mathcal J(s,z) D^2f(X(z))(\mathcal D_r X(z) u, \mathcal J(z,t)v )  dz \\
			&= \int^s_{t\vee r} \mathcal J(s,z) D^2f(X(z)) (\mathcal J(z,r) \sigma u, \mathcal J(z,t)v )  dz.
		\end{align*}
        For \(r>s\), this derivative is zero.
		Applying the product and chain rules to \(g(t)=\mathcal J(1,t)Df(X(t))\), we obtain
		\begin{align} \label{equ_Drg}
			(\mathcal D_r g(t) [u])v &= (\mathcal D_r \mathcal J(1,t) [u])(Df(X(t)) v)+  \mathcal J(1,t) D^2f(X(t)) (\mathcal D_r X(t)u,v)  \nn \\
			&=  \int^1_{t\vee r} \mathcal J(1,z) D^2f(X(z))\big(\mathcal J(z,r) \sigma u, \mathcal J(z,t)Df(X(t)) v \big)  dz \nn\\
			&\quad + \mathbf{1}_{r\leq t}  \mathcal J(1,t) D^2f(X(t)) (\mathcal J (t,r)\sigma u,v) \nn \\
			&=: G_{1,r}(t)(u,v) +G_{2,r}(t)(u,v).
		\end{align} 
		For the first term in \eqref{equ_Drg}, we have 
		\begin{align*}
			&\quad G_{1,r}(s)(u,v)- G_{1,r}(t)(u,v) \\
			&= \int^1_{s\vee r} \mathcal J(1,z)D^2f(X(z))\Big(\mathcal J(z,r) \sigma u,\big[ \mathcal J(z,s)Df(X(s)) -\mathcal J(z,t)Df(X(t))  \big] v \Big)  dz \\
			&\quad -\int_{t\vee r}^{s\vee r}
			\mathcal J(1,z) D^2f(X(z)) \big(
			\mathcal J(z,r)\sigma  u, \mathcal J(z,t)Df(X(t)) v\big)  dz \\
			&= \int^1_{s\vee r} \mathcal J(1,z) D^2f(X(z)) \Big(\mathcal J(z,r) \sigma u,\big[\mathcal J(z,s) - \mathcal J(z,t)  \big]Df(X(s))v \\
			&\qquad +\mathcal J(z,t) \big[Df(X(s))- Df(X(t)) \big] v \Big)  dz \\
			&\quad -\int_{t\vee r}^{s\vee r}
			\mathcal J(1,z)D^2f(X(z))
			\big(\mathcal J(z,r) \sigma u, \mathcal J(z,t)Df(X(t)) v\big)  dz. 
		\end{align*}
		The uniform bounds on \(\mathcal J\), \(Df\), and \(D^2f\), together with \(\eqref{equ_JTemporalHolder}\), therefore imply
		\begin{align}\label{equ_G1r}
			\| G_{1,r}(s)- G_{1,r}(t)\|_{\op} \leq C\big(|s-t| +|X(s)-X(t)|\big).
		\end{align}
		We now estimate the increments of $G_{2,r}$. For $r\in[0,t]$, 
		\begin{align*}
			&\quad G_{2,r}(s) (u,v)-G_{2,r}(t) (u,v) \\
			&= \mathcal J(1,s)D^2f(X(s))(\mathcal J(s,r)\sigma u,v) -\mathcal J(1,t)D^2f(X(t))(\mathcal J(t,r) \sigma u,v) \\
			&= \big( \mathcal J(1,s)-\mathcal J(1,t)\big) D^2f(X(s))(\mathcal J(s,r)\sigma u, v) \\
			&\quad+ \mathcal J(1,t)\big(  D^2f(X(s))- D^2f(X(t))\big) (\mathcal J(s,r)\sigma u,v)  \\
			&\quad +\mathcal J(1,t)  D^2f(X(t))\big((\mathcal J(s,r)- \mathcal J(t,r))\sigma u, v \big).  
		\end{align*} 
		For $r\in(t,s]$,
		\begin{align*}
			G_{2,r}(s) (u,v)-G_{2,r}(t) (u,v) =  \mathcal J(1,s)D^2f(X(s))(\mathcal J(s,r)\sigma u, v), \qquad r\in(t,s],
		\end{align*}
		whereas the difference is zero for $r\in(s,1]$.
		Similarly, the uniform boundedness and the temporal continuity of $\mathcal J$ lead to
		\begin{align}\label{equ_G2r}
			\| G_{2,r}(s)- G_{2,r}(t)\|_{\op} \leq C\big(|s-t| +|X(s)-X(t)|\big)+C\mathbf{1}_{\{r\in(t,s]\}}.
		\end{align}
		Combining \eqref{equ_G1r} and \eqref{equ_G2r}, and using the temporal continuity of $X$ \eqref{equ_XMomentcontinuity}, we have
		\begin{align*} 
			&\quad \Big(  \E \Big[ \big(\int^1_0 \big\| \mathcal D_rg(s)-\mathcal D_rg(t)\big\|_{\op}^{2} dr \big)^{\frac p2} \Big]\Big)^{\frac 1p} \nn\\
			&\leq  \Big(   \int^1_0  \big(\E\big[\big\| \mathcal D_rg(s)-\mathcal D_rg(t)\big\|_{\op}^{p}\big]\big)^{\frac 2p} dr  \Big)^{\frac 12} \nn\\ 
			&\leq  C_p\Big(   \int^1_0  \big(\E\big[\big\| G_{1,r}(s)-G_{1,r}(t)\big\|_{\op}^{p}+\big\| G_{2,r}(s)-G_{2,r}(t)\big\|_{\op}^{p}\big]\big)^{\frac 2p} dr  \Big)^{\frac 12} \nn\\ 
			&\leq C_p \Big(|s-t|^{\frac 12} +\big(\int^1_0 \mathbf{1}_{\{r\in(t,s]\}} dr\big)^{\frac 12}\Big) \leq C_p|s-t|^{\frac 12}.
		\end{align*}
		Here we use  the Minkowski integral inequality $\| \int^1_0 F(r,\cdot) dr \|_{L^{p/2}(\Omega)}\leq \int^1_0 \| F(r,\cdot)\|_{L^{ p/2}(\Omega)}dr$ for any $p\geq 2$ and any measurable function $F$ satisfying
		$\int_0^1 \|F(r,\cdot)\|_{L^{p/2}(\Omega)}\,\mathrm dr<\infty$, where $\|\cdot\|_{L^p(\Omega)}:= (\E [|\cdot|^p] )^{\frac 1p}$.

	\subsection{Proofs of the results in Section~\ref{sec_Comparison}}\label{subsec_proofinseccompare}
\begin{proof}[Proof of \eqref{equ_asO1}]
Let $h_i=N_i^{-1}$.
	 Applying 
		\eqref{equ_thetaEM} to the linear equation \eqref{equ_OU}, we have
		\begin{align*}
			X^{N_i,\theta}(1) = \Big( \frac{1+(1-\theta)\mu h_i }{1-\theta\mu h_i }\Big)^{N_i} x_0 + \sigma \frac{1}{1-\theta \mu h_i } \int^1_0 \Big( \frac{1+(1-\theta)\mu h_i }{1-\theta\mu h_i } \Big)^{N_i-1-[N_is]} dW(s).
		\end{align*}
		For $i\in\N^+$ large enough, the Taylor expansion gives
		$$
		\log \Big(  \frac{1+(1-\theta)\mu h_i }{1-\theta\mu h_i } \Big)
		=\mu h_i
		+\big(\theta-\frac12 \big)\mu^2 h_i^2
		+O(h_i^{3}),
		$$
		which leads to 
		\begin{align*}
		N_i	 \Big(  \Big( \frac{1+(1-\theta_1)\mu h_i }{1-\theta_1\mu h_i } \Big)^{N_i} -\Big(\frac{1+(1-\theta_2)\mu h_i }{1-\theta_2\mu h_i }\Big)^{N_i} \Big) = (\theta_1-\theta_2)\mu^2 e^\mu + O(h_i)
		\end{align*}
		and
			\begin{align*}
	&\quad  \frac{1}{1-\theta\mu h_i }	\Big( \frac{1+(1-\theta)\mu h_i }{1-\theta\mu h_i } \Big)^{N_i-1-[N_is]} \\
		&	=e^{\mu (1-s)}
		+e^{\mu (1-s)} h_i
		\Big[
		\mu\big(\{N_i s\}-\frac12 \big)
		+\big(\theta-\frac12 \big)\big(\mu+\mu^2(1-s)\big)
		\Big]
		+O(h_i^{2})
			\end{align*}
		for any $s\in[0,1]$. 
		Consequently, we have 
		\begin{align*}
			N_i (X^{N_i,\theta_1}(1)-X^{N_i,\theta_2}(1)) = C^{\theta_1,\theta_2}+ R_i,
		\end{align*}
		where \[
       R_i:= a_i +\int^1_0 b_i(s)dW(s),\qquad |a_i|+\sup_{s\in[0,1]} |b_i(s)| \leq Ch_i,
        \]
      with deterministic $a_i$ and $b_i$. The It\^o isometry therefore yields 
		 $\mathbb{E}[|R_i|^2] \leq C( N_i)^{-2}$. Since $\sum_{i}( N_i)^{-2}<\infty$, the Chebyshev inequality and Borel--Cantelli lemma imply that $R_i\to 0$ almost surely, which completes the proof.
\end{proof}

\begin{proof}[Proof of LILs for the $\theta$-Euler--Maruyama method applied to \eqref{eq:origin}] 
	 Using $Df(X(t)) \equiv A$, $\mathcal J(t,s)= e^{A(t-s)}$, and $\sigma = (\lambda, 0)^\top$, we derive
	 	\begin{align*}
	 		 \mathcal V_\sigma &= \frac{1}{12}  \int^1_0\mathcal J(1,r) Df(X(r))\sigma\sigma^\top  (Df(X(r)))^\top (\mathcal J(1,r))^\top  dr \\
	 		 &=\frac{1}{12}  \int^1_0 e^{A(1-r)} \begin{pmatrix}
	 		 	-\mu & -\gamma \\ 1 & 0
	 		 \end{pmatrix}  \begin{pmatrix}
	 		 \lambda \\ 0
	 		 \end{pmatrix}  \begin{pmatrix}
	 		\lambda & 0
	 		 \end{pmatrix} \begin{pmatrix}
	 		 -\mu & 1 \\ -\gamma & 0
	 		 \end{pmatrix} (e^{A(1-r)})^\top  dr \\ 
	 		 &= \frac{\lambda^2}{12}  \int^1_0 e^{At}  \begin{pmatrix}
	 		 	\mu^2 & -\mu \\ -\mu & 1
	 		 \end{pmatrix}   (e^{At})^\top  dt.
	 	\end{align*}
	 	For $\mu=2$ and $\gamma=1$, we write $A=-I_2 + (\begin{smallmatrix}
	 		-1 & -1 \\ 1 & 1
	 	\end{smallmatrix}) =: -I_2 +B$. Since $B^2=0$, we have
	 		\begin{align*}
	 	e^{At} = e^{-t}(I_2 +Bt) =e^{-t}  \begin{pmatrix}
	 	1-t & -t \\ t & 1+t
	 	\end{pmatrix},
	 		\end{align*}
	 		which leads to
	 		\begin{align*}
	 			 \mathcal V_\sigma &=  \frac{\lambda^2}{48} \begin{pmatrix}
	 			 	5-e^{-2} & -2 \\ -2 & 1-e^{-2}
	 			 \end{pmatrix}
	 		\end{align*} and $ \lambda_{\max}(\mathcal V_{\sigma} ) =  \frac{\lambda^2}{48}( 3-e^{-2}+2\sqrt{2})\approx 0.1186 \lambda^2$. 
	Let $X^{N_i,\theta}:=(P^{N_i,\theta}, Q^{N_i,\theta})^\top$  denote the numerical solution generated by the \(\theta\)-Euler–Maruyama method, and let \(e_p=(1,0)^\top\) and \(e_q=(0,1)^\top\). It follows from Remark~\ref{thm_thetamethod} that
		\begin{align*}
		\limsup_{i \to \infty} \frac{{N_i}(P(1)-P^{N_i,\theta}(1))  }{\sqrt{\log\log N_i}} 
		&=\sqrt{2 e_p^\top \mathcal V_{\sigma}  e_p}= \frac{\lambda}{2\sqrt{6}} \sqrt{5-e^{-2}},\quad\text{a.s.},
	\end{align*}
	\begin{align*}
		\limsup_{i \to \infty} \frac{{N_i}(Q(1)-Q^{N_i,\theta}(1))  }{\sqrt{\log\log N_i}} 
		&=\sqrt{2 e_q^\top \mathcal V_{\sigma} e_q}= \frac{\lambda}{2\sqrt{6}} \sqrt{1-e^{-2}},\quad\text{a.s.},
\end{align*}
and
		\begin{align*}
		\limsup_{i \to \infty} \frac{{N_i}| X(1)-X^{N_i,\theta}(1) | }{\sqrt{\log\log N_i}} = \sqrt{2\lambda_{\max}(\mathcal V_{\sigma} )}= \frac{\lambda}{2\sqrt{6}} \sqrt{3-e^{-2}+2\sqrt{2}} \approx 0.4870\lambda ,\quad\text{a.s.}
	\end{align*}
    The proof is finished.
\end{proof}

\begin{proof}[Proof of LILs for the BAOAB method applied to \eqref{eq:origin}] Let $h=N^{-1}$, $t_n=nh$, and $t_{n+\frac12}=t_n+h/2$.
The BAOAB method is defined by the composition $$(P^S_{n+1}, Q^S_{n+1}) = \Phi^B_{t_{n+1/2}, t_{n+1}}  \circ \Phi^A_{t_{n+1/2}, t_{n+1}} \circ \Phi^O_{t_n,t_{n+1}}  \circ \Phi^A_{t_n, t_{n+1/2}}\circ \Phi^B_{t_n, t_{n+1/2}} (P^S_{n}, Q^S_n),$$ 
where \(\Phi^A\), \(\Phi^B\), and \(\Phi^O\) are the flows of the free motion, force impulse, and Ornstein–Uhlenbeck subsystems, respectively:
\[
\begin{aligned}
A:\quad &dP=0, &&dQ=P\,dt,\\
B:\quad &dP=-\gamma Q\,dt, &&dQ=0,\\
O:\quad &dP=-\mu P\,dt+\lambda\,dW, &&dQ=0.
\end{aligned}
\]

Evaluating these flows successively gives
\begin{align}\label{equ_QBAOonestep}
		Q_{n+1}^S
		&= Q_n^S + \frac{h}{2} (1+ e^{-\mu h}) \big(P_n^S- \frac{h}{2} \gamma Q_n^S \big)+ \frac{\lambda h}{2} \int_{t_n}^{t_{n+1}} e^{-\mu(t_{n+1}-s)} dW(s), 
	\end{align}
	and
		\begin{align}\label{equ_PBAOonestep}
		P_{n+1}^S 
		&= e^{-\mu h} \big( P_n^S - \frac{h}{2} \gamma Q_n^S \big)+ \lambda \int_{t_n}^{t_{n+1}} e^{-\mu(t_{n+1}-s)} dW(s) -\frac{h}{2} \gamma Q_{n+1}^S  \nn\\
		&= \Big[e^{-\mu h} -\frac{\gamma h^2}{4}  (1+ e^{-\mu h})  \Big] P_n^S  -\frac{\gamma h}{2} (1+ e^{-\mu h})   (1- \frac{\gamma h^2}{4})  Q_n^S  \nn \\
		&\quad +\lambda (1- \frac{\gamma h^2}{4})  \int_{t_n}^{t_{n+1}} e^{-\mu(t_{n+1}-s)} dW(s). 
	\end{align}
Denote $X_n^S :=(P_n^S ,Q_n^S )^{\top}$.	Based on the one-step expression \eqref{equ_QBAOonestep}--\eqref{equ_PBAOonestep}, we define the continuous-time extension by  
	\begin{align} \label{equ_ContinuationBAOAB}
		X^{N,S}(t) & =M(t-t_n)X_n^S + \lambda b(t-t_n)\, \Xi^N_n(t), \qquad t\in[t_n,t_{n+1}],
	\end{align}
	where
	\begin{align*}
	M(\tau):=	\Bigg( \begin{matrix}
			e^{-\mu\tau}-\frac{\gamma\tau^2}{4}(1+e^{-\mu\tau})
			&
			-\dfrac{\gamma\tau}{2}(1+e^{-\mu\tau})
			(1-\dfrac{\gamma\tau^2}{4})
			\\
			\frac{\tau}{2}(1+e^{-\mu\tau})
			&
			1-\frac{\gamma\tau^2}{4}(1+e^{-\mu\tau})
		\end{matrix}\Bigg), 	\qquad b(\tau):= \begin{pmatrix} 
		1-\frac{\gamma\tau^2}{4} \\[5pt] \frac{\tau}{2}
		\end{pmatrix} 
	\end{align*}
	and   $\Xi^N_n(t): = \int_{t_n}^{t} e^{-\mu(t-s)} dW(s)$ for $t\in [t_n,t_{n+1}]$. 
	Define the error process $\mathcal E_{N,S}(t):= X(t)-X^{N,S}(t)$.
	Since  \eqref{eq:origin} is linear, its solution has the explicit representation 
	\begin{align} \label{equ_explicitLangevin}
		X(t) 
		=e^{At}X_0 
		+\lambda \int_{0}^t e^{A(t-s)} e_p  dW(s).
	\end{align}
	Combining  \eqref{equ_ContinuationBAOAB} and \eqref{equ_explicitLangevin},    we obtain, for $t\in[t_n,t_{n+1}]$,
	\begin{align*}
		&\quad \mathcal E_{N,S}(t) \\
        &=\Big( e^{A(t-t_n)}X(t_n)
		+\lambda\int_{t_n}^t e^{A(t-s)}e_p \, dW(s)\Big) -\Big( M(t-t_n) X_n^S +\lambda b(t-t_n)\, \Xi^N_n(t)\Big)  \\
		& = e^{A(t-t_n)}	\mathcal E_{N,S}(t_n) + \big( e^{A(t-t_n)}-M(t-t_n)\big)X_n^S+ \lambda \int_{t_n}^t\big(e^{A(t-s)} e_p-   b(t-t_n) e^{-\mu(t-s)}\big) dW(s).
	\end{align*}
	 By iteration, we have
	 	\begin{align*}
	 	\mathcal E_{N,S}(1) &=\sum_{n=0}^{N-1}
	 	e^{A(1-t_{n+1})}\big(e^{Ah}-M(h)\big)X_n^S\\
	 	&\quad +\lambda\sum_{n=0}^{N-1}
	 	\int_{t_n}^{t_{n+1}}
	 	e^{A(1-t_{n+1})}
	 	\big(
	 	e^{A(t_{n+1}-s)}e_p
	 	-b(h)e^{-\mu(t_{n+1}-s)}
	 	\big)dW(s) \\
	 	&=\sum_{n=0}^{N-1}
	 	e^{A(1-t_{n+1})}\big(e^{Ah}-M(h)\big)X_n^S\\
	 	&\quad +\lambda\sum_{n=0}^{N-1}
	 	\int_{t_n}^{t_{n+1}}
	 	e^{A(1-s)}\big(
	 	e_p
	 	-e^{-(A+\mu I_2)(t_{n+1}-s)}b(h)
	 	\big)dW(s). 
	 \end{align*} 
	 It follows from $(A+\mu I_2) e_p =e_q$ and $b(h) = (1-\frac{\gamma h^2}{4})e_p +\frac h2 e_q$ that
	 \begin{align}\label{equ_expanStochasticBAOAB}
	 &\quad 	e_p
	 	-e^{-(A+\mu I_2)(t_{n+1}-s)}b(h) \nn\\
	 	& =e_p - (I_2-(A+\mu I_2)(t_{n+1}-s)+\mathcal {O} ( (t_{n+1}-s)^2)) \Big( (1-\frac{\gamma h^2}{4})e_p +\frac h2 e_q\Big) \nn \\
	 	& = \Big(t_{n+1}-s -\frac h2 \Big) e_q +\mathcal{O}(h^2),
	 \end{align}
     uniformly in $n$ and $s\in[t_n,t_{n+1}]$. 
	 Combining \eqref{equ_expanStochasticBAOAB} with the fact that $\| e^{Ah}-M(h)\|_{\op} \leq Ch^3$ and 
     using $t_{n+1}-s-\frac h2
=h\left(\frac12-\{Ns\}\right)$
for almost every $s\in[t_n,t_{n+1}]$,
we therefore obtain $N{\mathcal E}_{N,S}(1)=Y_N^S+R_N^S$,
where $(\E[|R_N^S|^2])^{\frac 12}\leq CN^{-1}$ and \begin{align*}
	 	Y^{S}_N:=\lambda \int^1_0 \mathcal J(1,r) e_q \big(\frac 12-\{ Nr\}\big)dW(r)
	 \end{align*}
with $\mathcal J(1,r)=e^{A(1-r)}$.
Since \(\sum_i^\infty (N_i)^{-2}<\infty\), the Chebyshev inequality and the Borel–Cantelli lemma yield \(R_{N_i}^S\to0\) almost surely. Thus, the LIL asymptotics of the terminal error are determined by \(Y_{N_i}^S\).

Applying the argument used in the proof of \Cref{prop_LILofY} to \(Y_N^S\), we obtain the corresponding LILs with covariance matrix	
	 \begin{align*}
	 	\mathcal V^S_{\sigma}= \frac{\lambda^2 }{12} \int^1_0 e^{A(1-r)} e_q e_q^\top (e^{A(1-r)})^\top dr =   \frac{\lambda^2 }{12} \int^1_0 e^{Ar} \begin{pmatrix}
	 		0 & 0 \\ 0 & 1
	 	\end{pmatrix}  (e^{Ar})^\top dr .
	 \end{align*}

	 
For \(\mu=2\) and \(\gamma=1\), evaluating the integral gives
\begin{align*}
	 	\mathcal V^S_\sigma &=  \frac{\lambda^2}{48} \begin{pmatrix}
	 		1-5e^{-2} & -2 +8e^{-2}\\ -2+8e^{-2} & 5-13e^{-2}
	 	\end{pmatrix},
	 \end{align*} 
     whose largest eigenvalue is $ \lambda_{\max}(\mathcal V^S_{\sigma} ) =  \frac{\lambda^2}{48}( 3-9e^{-2}+2\sqrt{2-12e^{-2}+20e^{-4}}) \approx 0.073 \lambda^2$. The LIL constants for the normalized momentum error, position error, and Euclidean error norm then follow as in the proof for the \(\theta\)-Euler–Maruyama method. The proof is complete.
    \end{proof}

	\bibliographystyle{plain}
	\bibliography{reference.bib}

@book {Milstein2004book,
	AUTHOR = {Milstein, G. N. and Tretyakov, M. V.},
	TITLE = {Stochastic {N}umerics for {M}athematical {P}hysics},
	SERIES = {Scientific Computation},
	PUBLISHER = {Springer-Verlag, Berlin},
	YEAR = {2004},
	PAGES = {xx+594},
	ISBN = {3-540-21110-1},
	MRCLASS = {60-02 (35R60 37M15 60H10 60H15 60H35 65C30)},
	MRNUMBER = {2069903},
	MRREVIEWER = {Peter\ E.\ Kloeden},
	DOI = {10.1007/978-3-662-10063-9},
	URL = {https://doi.org/10.1007/978-3-662-10063-9},
}

@incollection {KurtzProtter1991SA,
	AUTHOR = {Kurtz, T. G. and Protter, P.},
	TITLE = {Wong-{Z}akai corrections, random evolutions, and simulation
	schemes for {SDE}s},
	BOOKTITLE = {Stochastic analysis},
	PAGES = {331--346},
	PUBLISHER = {Academic Press, Boston, MA},
	YEAR = {1991},
	ISBN = {0-12-481005-5},
	MRCLASS = {60H10 (60F17)},
	MRNUMBER = {1119837},
	MRREVIEWER = {Leszek\ S\l omi\'nski},
}

@book {Nualart2018book,
	AUTHOR = {Nualart, D. and Nualart, E.},
	TITLE = {Introduction to {M}alliavin calculus},
	SERIES = {Institute of Mathematical Statistics Textbooks},
	VOLUME = {9},
	PUBLISHER = {Cambridge University Press, Cambridge},
	YEAR = {2018},
	PAGES = {xii+236},
	ISBN = {978-1-107-61198-6; 978-1-107-03912-4},
	MRCLASS = {60-01 (60H07 60J65 60J75)},
	MRNUMBER = {3838464},
	MRREVIEWER = {Maria\ Gordina},
	DOI = {10.1017/9781139856485},
	URL = {https://doi.org/10.1017/9781139856485},
}

@article {ConditionalKSthm2022,
	AUTHOR = {Lo, C. H. and Menshikov, M. V. and Wade, A. R.},
	TITLE = {Cutpoints of non-homogeneous random walks},
	JOURNAL = {ALEA Lat. Am. J. Probab. Math. Stat.},
	FJOURNAL = {ALEA. Latin American Journal of Probability and Mathematical
	Statistics},
	VOLUME = {19},
	YEAR = {2022},
	NUMBER = {1},
	PAGES = {493--510},
	ISSN = {1980-0436},
	MRCLASS = {60J05 (60G50 60J10)},
	MRNUMBER = {4394306},
	MRREVIEWER = {Barbara\ Martinucci},
	DOI = {10.30757/alea.v19-19},
	URL = {https://doi.org/10.30757/alea.v19-19},
}

@article {KochenStone1964,
	AUTHOR = {Kochen, S. and Stone, C.},
	TITLE = {A note on the {B}orel-{C}antelli lemma},
	JOURNAL = {Illinois J. Math.},
	FJOURNAL = {Illinois Journal of Mathematics},
	VOLUME = {8},
	YEAR = {1964},
	PAGES = {248--251},
	ISSN = {0019-2082},
	MRCLASS = {60.30},
	MRNUMBER = {161355},
	MRREVIEWER = {K.\ L.\ Chung},
	URL = {http://projecteuclid.org/euclid.ijm/1256059668},
}

@article {Bichteler1981AOP,
	AUTHOR = {Bichteler, K.},
	TITLE = {Stochastic integration and {$L\sp{p}$}-theory of
	semimartingales},
	JOURNAL = {Ann. Probab.},
	FJOURNAL = {The Annals of Probability},
	VOLUME = {9},
	YEAR = {1981},
	NUMBER = {1},
	PAGES = {49--89},
	ISSN = {0091-1798,2168-894X},
	MRCLASS = {60H05 (60G48 60H20)},
	MRNUMBER = {606798},
	MRREVIEWER = {Philip\ Protter},
	URL =
	{http://links.jstor.org/sici?sici=0091-1798(198102)9:1<49:SIAOS>2.0.CO;2-X&origin=MSN},
}

@article {ProtterSPA2020,
	AUTHOR = {Protter, P. and Qiu, L. and Martin, J. S.},
	TITLE = {Asymptotic error distribution for the {E}uler scheme with
	locally {L}ipschitz coefficients},
	JOURNAL = {Stochastic Process. Appl.},
	FJOURNAL = {Stochastic Processes and their Applications},
	VOLUME = {130},
	YEAR = {2020},
	NUMBER = {4},
	PAGES = {2296--2311},
	ISSN = {0304-4149,1879-209X},
	MRCLASS = {60H35 (60F05 60H10)},
	MRNUMBER = {4074703},
	MRREVIEWER = {C\'onall\ Kelly},
	DOI = {10.1016/j.spa.2019.07.003},
	URL = {https://doi.org/10.1016/j.spa.2019.07.003},
}

@article {JacodProtter1998AOP,
	AUTHOR = {Jacod, J. and Protter, P.},
	TITLE = {Asymptotic error distributions for the {E}uler method for
	stochastic differential equations},
	JOURNAL = {Ann. Probab.},
	FJOURNAL = {The Annals of Probability},
	VOLUME = {26},
	YEAR = {1998},
	NUMBER = {1},
	PAGES = {267--307},
	ISSN = {0091-1798,2168-894X},
	MRCLASS = {60H10 (60G44 65U05)},
	MRNUMBER = {1617049},
	MRREVIEWER = {Denis\ Talay},
	DOI = {10.1214/aop/1022855419},
	URL = {https://doi.org/10.1214/aop/1022855419},
}

@book{schilling2014brownian,
	AUTHOR = {Schilling, R. L. and Partzsch, L.},
	TITLE = {{B}rownian {M}otion: an {I}ntroduction to {S}tochastic {P}rocesses},
	SERIES = {De Gruyter Graduate},
	EDITION = {Second},
	PUBLISHER = {De Gruyter, Berlin},
	YEAR = {2014},
	PAGES = {xvi+408},
	ISBN = {978-3-11-030729-0; 978-3-11-030730-6},
	MRCLASS = {60J65 (60G05 60G15 60G17 60G44 60H05 60H10 60J25)},
	MRNUMBER = {3234570},
	MRREVIEWER = {Tam\'{a}s\ Szabados},
	URL = {https://doi.org/10.1515/9783110307306},
}

@book {Leimkuhler2015book,
	AUTHOR = {Leimkuhler, B. and Matthews, C.},
	TITLE = {Molecular dynamics},
	SERIES = {Interdisciplinary Applied Mathematics},
	VOLUME = {39},
	NOTE = {With deterministic and stochastic numerical methods},
	PUBLISHER = {Springer, Cham},
	YEAR = {2015},
	PAGES = {xxii+443},
	ISBN = {978-3-319-16374-1; 978-3-319-16375-8},
	MRCLASS = {82Bxx (35Q82 65-02 70-02 82Cxx)},
	MRNUMBER = {3362507},
}

@article {Gyongy1998PA,
	AUTHOR = {Gy\"ongy, I.},
	TITLE = {A note on {E}uler's approximations},
	JOURNAL = {Potential Anal.},
	FJOURNAL = {Potential Analysis. An International Journal Devoted to the
	Interactions between Potential Theory, Probability Theory,
	Geometry and Functional Analysis},
	VOLUME = {8},
	YEAR = {1998},
	NUMBER = {3},
	PAGES = {205--216},
	ISSN = {0926-2601,1572-929X},
	MRCLASS = {60H10},
	MRNUMBER = {1625576},
	MRREVIEWER = {Shigetoku\ Kawabata},
	DOI = {10.1023/A:1008605221617},
	URL = {https://doi.org/10.1023/A:1008605221617},
}

@article {Shardlow2016BIT,
	AUTHOR = {Shardlow, T. and Taylor, P.},
	TITLE = {On the pathwise approximation of stochastic differential
	equations},
	JOURNAL = {BIT},
	FJOURNAL = {BIT. Numerical Mathematics},
	VOLUME = {56},
	YEAR = {2016},
	NUMBER = {3},
	PAGES = {1101--1129},
	ISSN = {0006-3835,1572-9125},
	MRCLASS = {65C30},
	MRNUMBER = {3540473},
	MRREVIEWER = {John\ Masson\ Noble},
	DOI = {10.1007/s10543-015-0597-2},
	URL = {https://doi.org/10.1007/s10543-015-0597-2},
}

@article {Kloeden2007LMS,
	AUTHOR = {Kloeden, P. E. and Neuenkirch, A.},
	TITLE = {The pathwise convergence of approximation schemes for
	stochastic differential equations},
	JOURNAL = {LMS J. Comput. Math.},
	FJOURNAL = {LMS Journal of Computation and Mathematics},
	VOLUME = {10},
	YEAR = {2007},
	PAGES = {235--253},
	ISSN = {1461-1570},
	MRCLASS = {60H35 (60H10 65C30)},
	MRNUMBER = {2320830},
	MRREVIEWER = {Xue\ Rong\ Mao},
	DOI = {10.1112/S1461157000001388},
	URL = {https://doi.org/10.1112/S1461157000001388},
}

@article{Chenchen2024,
      title={Superiority of stochastic symplectic methods via the law of iterated logarithm}, 
      author={C. Chen and X. Chen and T. Dang and J. Hong},
      year={2024},
      journal = {\rm{arXiv:2404.14842}}
}

@article{Chenchen2025,
      title={The law of iterated logarithm for numerical approximation of time-homogeneous Markov process}, 
      author={C. Chen and X. Chen and J. Hong},
      year={2025},
      journal = {\rm{arXiv:2511.05217}}
}

@article {Huang2025,
    AUTHOR = {Chen, C. and Hong, J. and Huang, C.},
     TITLE = {Stochastic modified equations for symplectic methods applied
              to rough {H}amiltonian systems},
   JOURNAL = {IMA J. Numer. Anal.},
  FJOURNAL = {IMA Journal of Numerical Analysis},
    VOLUME = {45},
      YEAR = {2025},
    NUMBER = {2},
     PAGES = {894--933},
      ISSN = {0272-4979,1464-3642},
   MRCLASS = {65C30 (65P10)},
  MRNUMBER = {4886514},
       DOI = {10.1093/imanum/drae019},
       URL = {https://doi.org/10.1093/imanum/drae019},
}
	
\end{document}